%% file: contmm_energyanalysis.tex
\documentclass[12pt]{article}

\usepackage[longpaper]{preprint-layout}
\usepackage{preprint-notation}
\usepackage{preprint-tools}
\usepackage{hyperref}

\usepackage{dsfont}
\usepackage{import}
\usepackage{wrapfig}

\title{Space-Time CutFEM on Overlapping Meshes: Simple Continuous Mesh Motion}
\author{Mats G. Larson\footnote{\addressumu}, Anders Logg\footnote{\addressch}, Carl Lundholm\footnote{\addressumu}}
\date{\today}

\numberwithin{equation}{section}

\newtheorem{theorem}{Theorem}[section] 
\newtheorem{lemma}{Lemma}[section]
\newtheorem{corollary}{Corollary}[section]

\newcommand{\bnab}{\bar{\nabla}}

\newcommand{\bx}{\bar{x}}
\newcommand{\bs}{\bar{s}}
\newcommand{\bn}{\bar{n}}
\newcommand{\bm}{\bar{\mu}}
\newcommand{\abm}{|\bar{\mu}|}
\newcommand{\ud}{\,\mathrm{d}}
\newcommand{\nab}{\nabla}

\newcommand{\lap}{\Delta}
\newcommand{\sgn}{\text{sgn}}
\newcommand{\euler}{\text{e}}
\newcommand{\norma}[1]{\left|\mkern-1.5mu\left|\mkern-1.5mu\left| #1 \right|\mkern-1.5mu\right|\mkern-1.5mu\right|_{A_{h,t}}} 
\newcommand{\normb}[1]{\left|\mkern-1.5mu\left|\mkern-1.5mu\left| #1 \right|\mkern-1.5mu\right|\mkern-1.5mu\right|_{B_h}}
\newcommand{\normx}[1]{\left|\mkern-1.5mu\left|\mkern-1.5mu\left| #1 \right|\mkern-1.5mu\right|\mkern-1.5mu\right|_{X}}
\newcommand{\normym}[1]{\left|\mkern-1.5mu\left|\mkern-1.5mu\left| #1 \right|\mkern-1.5mu\right|\mkern-1.5mu\right|_{Y_-}}
\newcommand{\normyp}[1]{\left|\mkern-1.5mu\left|\mkern-1.5mu\left| #1 \right|\mkern-1.5mu\right|\mkern-1.5mu\right|_{Y_+}}
 
\newcommand{\bGn}{\bar{\Gamma}_n}

\newcommand{\Om}[1]{{\Omega_#1}}
\newcommand{\Omt}[1]{{\Omega_#1(t)}}
   
\newcommand{\interp}[2]{\bar{I}_{#1}^{#2}}
\newcommand{\interph}{\bar{I}_h}

\usetikzlibrary{arrows}
\usetikzlibrary{shapes}

\begin{document}

\maketitle

\begin{abstract}
We present a cut finite element method for the heat equation on two overlapping meshes: a stationary background mesh and an overlapping mesh that moves around inside/``on top'' of it. Here the overlapping mesh is prescribed a simple continuous motion, meaning that its location as a function of time is \emph{continuous} and \emph{piecewise linear}. For the discrete function space, we use continuous Galerkin in space and discontinuous Galerkin in time, with the addition of a discontinuity on the boundary between the two meshes. The finite element formulation is based on Nitsche's method and also includes an integral term over the space-time boundary between the two meshes that mimics the standard discontinuous Galerkin time-jump term. The simple continuous mesh motion results in a space-time discretization for which standard analysis methodologies either fail or are unsuitable. We therefore employ what seems to be a relatively new energy analysis framework that is general and robust enough to be applicable to the current setting. The energy analysis consists of a stability estimate that is slightly stronger than the standard basic one and an a priori error estimate that is of optimal order with respect to both time step and mesh size. We also present numerical results for a problem in one spatial dimension that verify the analytic error convergence orders. \\
\end{abstract}

\vspace{1cm}


\section{Introduction}\label{sec:introduction}

\paragraph{Issue - Cost of mesh generation:} Generating computational meshes for numerically solving differential equations can be a computationally costly procedure. In practical applications the mesh generation can often represent a substantial amount of the total computation time. This is especially true for problems where the solution domain changes during the solve process, e.g., evolving geometry and shape optimization. With standard methods the mesh then has to be constantly checked for degeneracy and updated if needed, meaning a persisting meshing cost for the entire solve process. 

\paragraph{Remedy - CutFEM:} Cut finite element methods (CutFEMs) provide a way of decoupling the computational mesh from the problem geometry. This means that the same discretization can be used for a changing solution domain. CutFEMs can thus make remeshing redundant for problems with changing geometry but also for other applications involving meshing such as adaptive mesh refinement. The cost of CutFEMs is treating the mesh cells that are arbitrarily cut by the independent problem geometry.

\begin{wrapfigure}{r}{0.5\textwidth}
\centering
\includegraphics[width=0.5\textwidth]{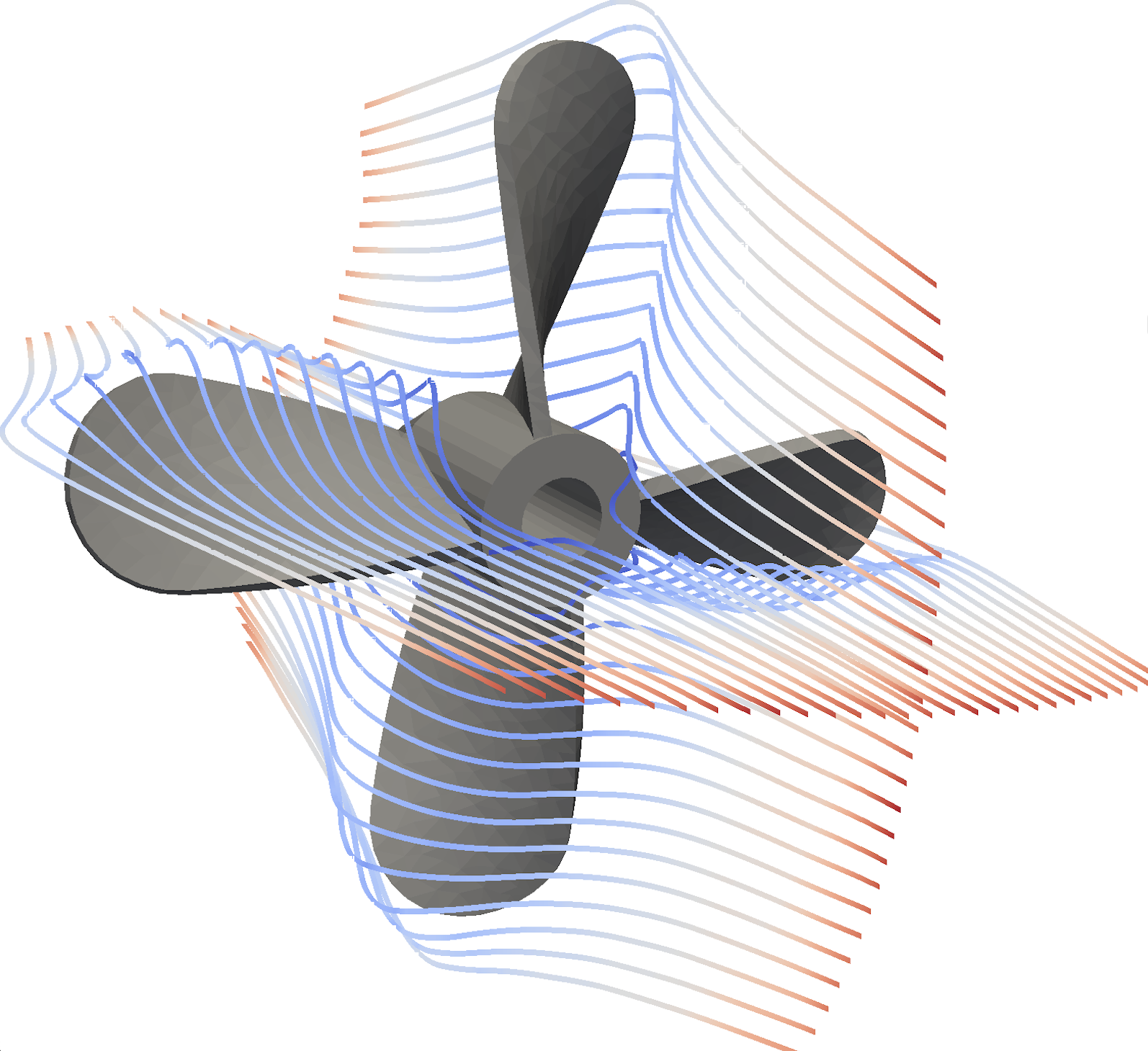}
\caption{Computed streamlines around a propeller. \href{https://link.springer.com/article/10.1186/s40323-015-0043-7}{Image} by \href{https://anders.logg.org/}{Anders Logg} is licensed under \href{https://creativecommons.org/licenses/by/4.0/}{CC BY 4.0}.}
\label{fig_propeller_flow}
\end{wrapfigure}

\begin{wrapfigure}{r}{0.5\textwidth}
\centering
\includegraphics[width=0.5\textwidth]{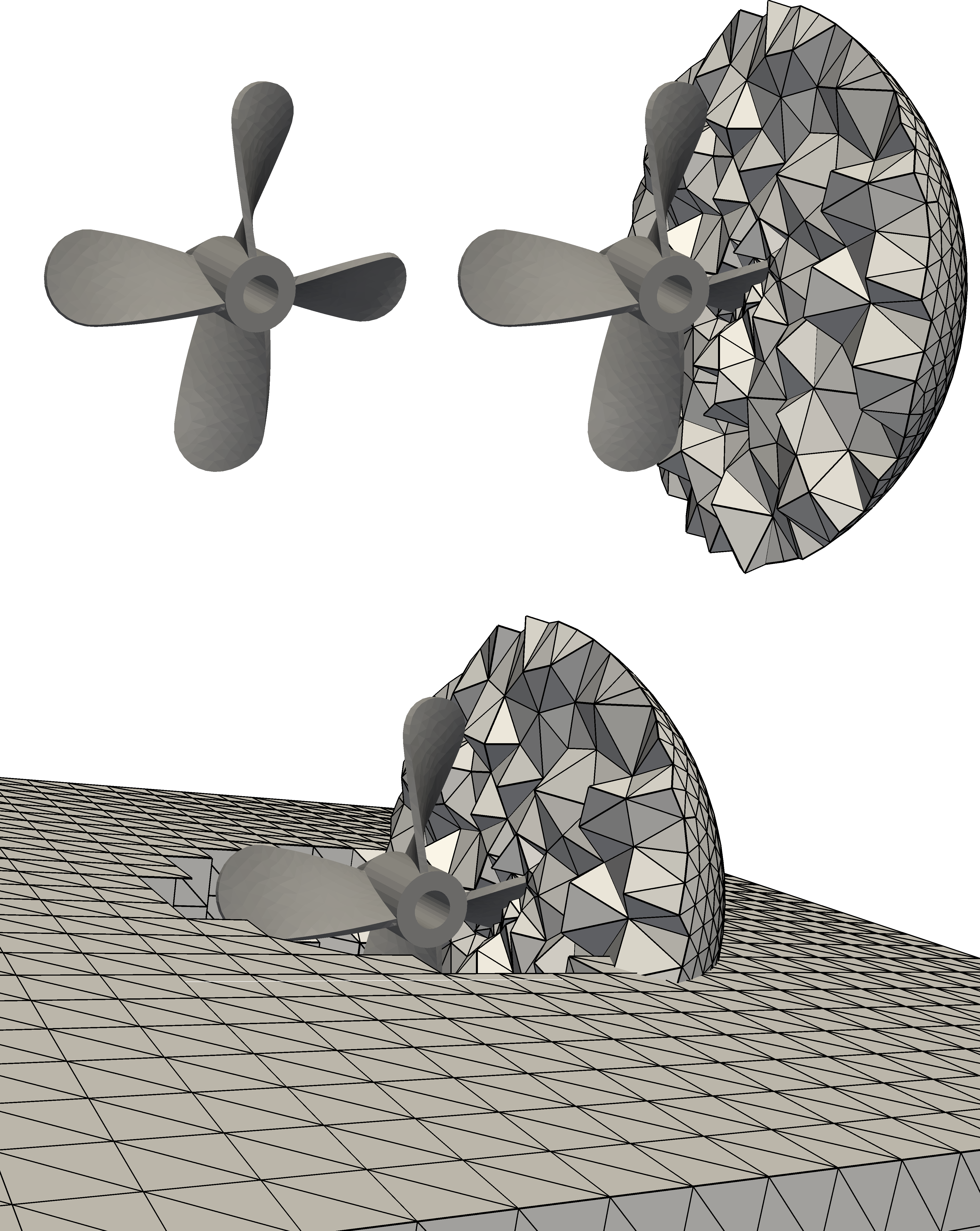}
\caption{Overlapping meshes for a problem with a rotating propeller. \href{https://link.springer.com/article/10.1186/s40323-015-0043-7}{Image} by \href{https://anders.logg.org/}{Anders Logg} is licensed under \href{https://creativecommons.org/licenses/by/4.0/}{CC BY 4.0}.}
\label{fig_propeller_om}
\end{wrapfigure}

\paragraph{CutFEM on overlapping meshes:} A common type of problem with changing geometry is one where there is a moving object in the solution domain, e.g., see Figure~\ref{fig_propeller_flow}. A straightforward CutFEM-approach to this problem would be to consider CutFEM for the interface problem, i.e., to use a background mesh of the empty solution domain together with an interface that represents the object. However, a more advantageous and sophisticated approach is to consider CutFEM on \emph{overlapping meshes}, meaning two or more meshes ordered in a mesh hierarchy. This is also called composite grids/meshes and multimesh in the literature. The idea is to use a background mesh of the empty solution domain, just as for the interface problem, but instead to encapsulate the object in a second mesh. The mesh containing the object is then placed ``on top'' of the background mesh, creating a mesh hierarchy. The motion of the object will thus also cause its encapsulating mesh to move. There are some advantages of using a second overlapping mesh instead of an interface. Firstly, an overlapping mesh can incorporate boundary layers close to the object. Something an interface cannot. Secondly, the total number of degrees of freedom (DOFs) of the resulting linear system may be reduced. This is so since for CutFEM for the interface problem this number can be twice the number of DOFs of the background mesh or more, whereas for CutFEM on overlapping meshes it will be the number of DOFs of the background mesh plus the number of DOFs of the second mesh. Thirdly, if the object has a complicated geometry, representing it with an interface can lead to tricky cut situations and thus a higher computational cost. By instead using an object-encapsulating mesh with a simply-shaped exterior boundary, the cut situations can be made less tricky, see Figure~\ref{fig_propeller_om}. A way to further sophisticate this is to allow the moving object to deform the interior of the overlapping mesh while initially keeping its exterior boundary fixed. Only when the deformations have become too large is the overlapping mesh ``snapped'' into place to avoid degeneracy. Such a snapping feature provides a choice between computing cut situations or computing deformations, thus allowing the cheapest option for the situation at hand to be chosen. A drawback of using a second overlapping mesh instead of an interface is that overlapping meshes require collision computations between the cells of the meshes, something that can be computationally expensive.

CutFEM on overlapping meshes can also be used as an alternative to adaptive mesh refinement by keeping a smaller finer mesh in regions requiring higher accuracy. Yet another application example is to use a composition of simpler structured meshes to represent a complicated domain. 

\paragraph{Literary background:} Over the past two decades, a theoretical foundation for the formulation of
stabilized CutFEM has been developed by extending the ideas
of Nitsche, presented in~\cite{Nitsche:1971aa}, to a general weak formulation of the
interface conditions, thereby removing the need for domain-fitted
meshes. The foundations of CutFEM were presented
in~\cite{Hansbo:2002aa} and then extended to overlapping meshes in~\cite{Hansbo:2003aa}. The CutFEM methodology has since been developed and applied to a number of important multiphysics problems.
See for example~\cite{Burman:2007ab,Burman:2007aa,Becker:2009aa,Massing2014a}. For overlapping meshes in particular, see for example~\cite{Massing2014,Johansson:2015aa, Dokken2019, Johansson2019}. So far, only CutFEM for \emph{stationary} problems on overlapping meshes have been developed and analysed to a satisfactory degree, thus leaving analogous work for \emph{time-dependent} problems to be desired.

\paragraph{This work:} The work presented here is intended to be an initial part of developing and analysing CutFEMs for time-dependent problems on overlapping meshes. We consider a CutFEM for the heat equation on two overlapping meshes: one stationary background mesh and one moving overlapping mesh with no object. Depending on how the mesh motion is represented discretely, quite different space-time discretizations may arise, allowing for different types of analyses to be applied. Generally the mesh motion may either be continuous or discontinuous. We have considered the simplest case of both of these two types, which we refer to as \emph{simple continuous} and \emph{simple discontinuous} mesh motion. Simple continuous mesh motion means that the location of the overlapping mesh as a function of time is \emph{continuous} and \emph{piecewise linear}, and simple discontinuous mesh motion means that it is \emph{discontinuous} and \emph{piecewise constant}. The latter is studied in other work and the former in this.

\paragraph{Analytic novelty:} The simple continuous mesh motion results in a space-time discretization with skewed space-time nodal trajectories and cut prismatic space-time cells. This discretization lacks a slabwise product structure between space and time. Standard analysis methodology relying on such a structure therefore either fail or require too restrictive assumptions here. The reason for this is that standard analysis methodology typically use spatial operators that map to the momentaneous finite element space, such as the discrete Laplacian and the solution operator used to define the $H^{-1}$-norm on $L^2$. If the spatial discretization changes within slabs these operators get an intrinsic time dependence that standard methodologies fail to incorporate. We therefore employ what seems to be a relatively new analysis framework for parabolic problems that avoids the use of operators of the aforementioned type which thus makes it general and robust enough to be applicable to the current discretization. It seems that the core components of this analysis framework have been discovered independently by us and Cangiani, Dong, and Georgoulis in \cite{Cangiani2017}. The analysis is of an energy type, where space-time energy norms are used to derive and obtain a stability estimate that is slightly stronger than the standard basic one and an a priori error estimate that is of optimal order with respect to both time step and mesh size. The main steps of this new energy analysis are:
\begin{enumerate}
\setcounter{enumi}{-1}
\item Handling of the time derivative: This is the initial step that characterizes and sets the course for the whole analysis. Instead of the $H^{-1}$-norm, the $L^2$-norm scaled with the time step is used to include the time-derivative term in a space-time energy norm.  
\item Analytic preliminaries: A ``perturbed coercivity'' is proved which is used to show an inf-sup condition. These results become slightly different compared with corresponding standard ones due to the new handling of the time derivative.
\item Stability analysis: The ``perturbed coercivity'' is used to derive a stability estimate that is somewhat stronger than the standard basic one obtained by testing with the discrete solution.
\item Error analysis: Just as in a standard energy analysis, a Cea's lemma type argument is followed by using the inf-sup condition, Galerkin orthogonality, and continuity. A difference here is that the continuity comes with a twist, namely temporal integration by parts, which is needed because of the slightly different inf-sup condition. Finally, together with an interpolation estimate, an optimal order a priori error estimate may be proved. 
\end{enumerate}

\paragraph{Paper overview:} In Section 2, the model problem is formulated. In Section 3, the CutFEM is presented. In Section 4, tools for the analysis are presented. In Section 5, we present and prove a stability estimate for the discrete solution. In Section 6, we present and prove an optimal order a priori error estimate. In Section 7, we present numerical results for a problem in one spatial dimension that verify the analytic convergence orders. In the appendix we present technical estimates and interpolation results used in the analysis.

\section{Problem} \label{sec:probform}

For $d = 1, 2$, or $3$, let $\Om{0} \subset \mathbb{R}^d$ be a bounded convex domain with polygonal boundary $\partial\Om{0}$. Let $T > 0$ be a given final time. Let $G \subset \Om{0} \subset \mathbb{R}^d$ be another bounded domain with polygonal boundary $\partial G$. We let the \emph{location} of $G$ be time-dependent by prescribing for $G$ a velocity $\mu : [0, T] \to \mathbb{R}^d$. We point out that the size and shape of $G$ remain the same for all times. Using $\Om{0}$ and $G$, we define the following two domains:
\begin{wrapfigure}{r}{0.5\textwidth}
\centering
\def\svgwidth{0.5\textwidth}
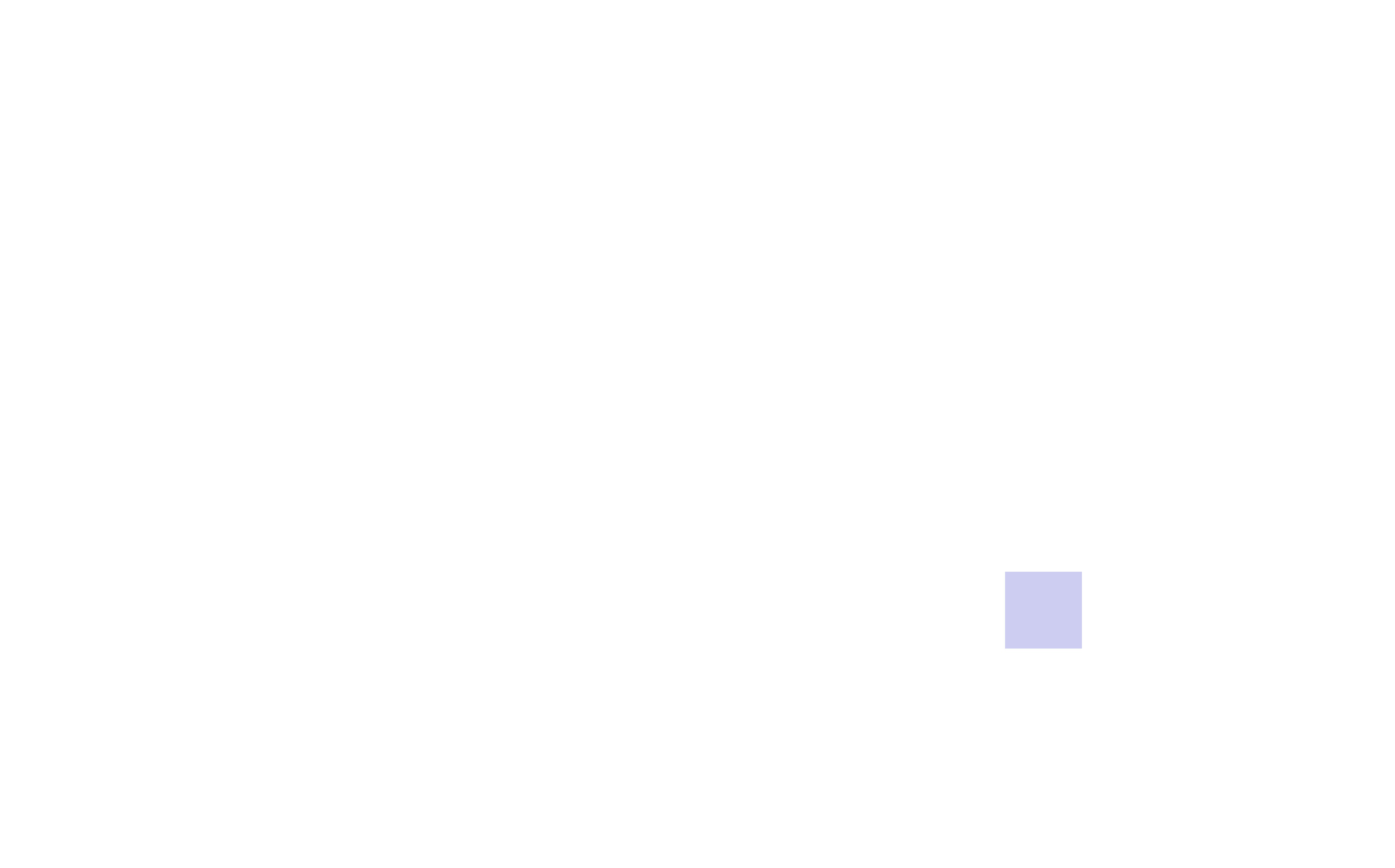
\caption{Partition of $\Om{0}$ into $\Om{1}$ (blue) and $\Om{2}$ (red) for $d = 2$ with $G$ moving with velocity $\mu$.}
\label{fig_probdom}
\end{wrapfigure}
\begin{align}
\Om{1} &:= \Om{0} \setminus (G \cup \partial G) \label{def:Om1} \\
\Om{2} &:= \Om{0} \cap G \label{def:Om2}
\end{align}
with boundaries $\partial\Om{1}$ and $\partial\Om{2}$, respectively. Let their common boundary be
\begin{equation}
\Gamma := \partial\Om{1} \cap \partial\Om{2} \label{def:Gamma}
\end{equation}
For $t \in [0, T]$, we have the partition
\begin{equation}
\Om{0} = \Omt{1} \cup \Gamma(t) \cup \Omt{2} \label{partitionOm0}
\end{equation}
See Figure~\ref{fig_probdom} for an illustration. We consider the heat equation in $\Om{0} \times (0, T]$ with source $f \in L^2((0, T], \Om{0})$, homogeneous Dirichlet boundary conditions, and initial data $u_0 \in H^2(\Om{0}) \cap H_0^1(\Om{0})$:
\begin{equation} 
\left\{
\begin{split}
\dot{u} - \lap{u}  & = f && \text{in} \ \Om{0} \times (0, T] \\
u & = 0 && \text{on} \ \partial\Om{0} \times (0, T] \\
u & = u_0 && \text{in} \ \Om{0} \times \{0\} 
\end{split}
\right. \label{heateq}
\end{equation}

\section{Method} \label{sec:feform}

\subsection{Preliminaries}

Let $\mathcal{T}_0$ and $\mathcal{T}_G$ be quasi-uniform simplicial meshes of $\Om{0}$ and $G$, respectively. We denote by $h_K$ the diameter of a simplex $K$. We partition the time interval $(0, T]$ quasi-uniformly into $N$ subintervals $I_n = (t_{n-1}, t_n]$ of length $k_n = t_n - t_{n-1}$, where $0 = t_0 < t_1 < \ \dots \ < t_N = T$ and $n = 1, \dots , N$. We assume the following space-time quasi-uniformity: For $h = \max_{K \in \mathcal{T}_0 \cup \mathcal{T}_G}\{h_K\}$, and $k = \max_{1 \leq n \leq N}\{k_n\}$,
\begin{equation}
h^2 \lesssim k_{\min} \quad k \lesssim h_{\min}
\label{quasiuniformity_st}
\end{equation}
where $k_{\min} = \min_{1 \leq n \leq N}\{k_n\}$, and $h_{\min} = \min_{K \in \mathcal{T}_0 \cup \mathcal{T}_G}\{h_K\}$. We next define the following slabwise space-time domains:
\begin{align}
S_{{0,n}} &:= \Om{0} \times I_n \label{def:S0n} \\
S_{{i,n}} &:= \{(x, t) \in S_{{0,n}} : x \in \Omt{i}\} \label{def:Sin} \\
\bGn &:= \{(s, t) \in S_{{0,n}} : s \in \Gamma(t)\} \label{def:bGn}
\end{align}
\begin{figure}[h]
\centering
\def\svgwidth{0.49\textwidth}
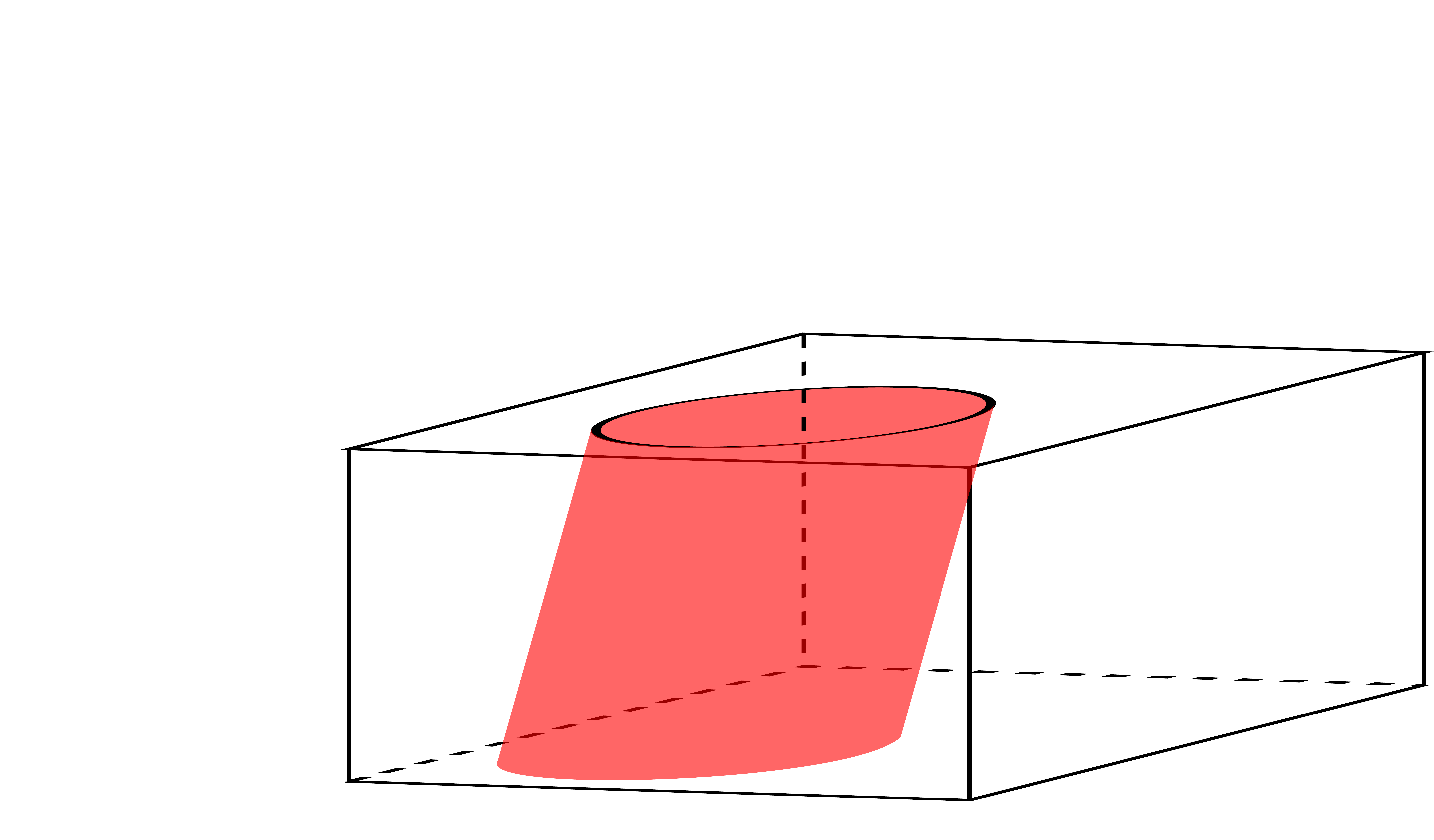
\def\svgwidth{0.49\textwidth}
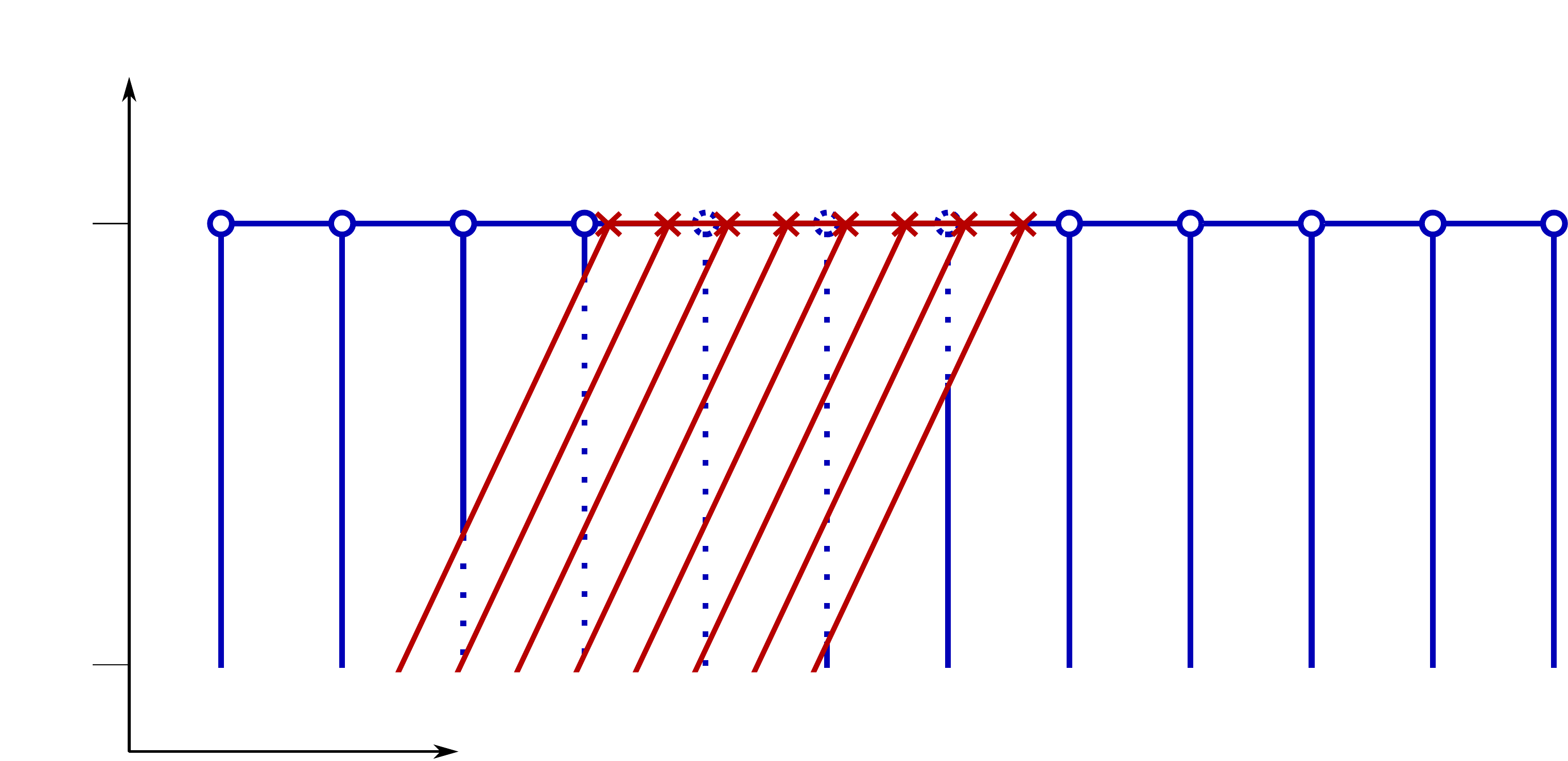
\caption{\textbf{Left:} Space-time slabs with simple continuous mesh motion. \textbf{Right:} Space-time discretization for $S_{{0,n}}$ for $d = 1$ when $\mu > 0$. At time $t = t_n$, the nodes of the blue background mesh $\mathcal{T}_0$ are marked with circles and the nodes of the red moving mesh $\mathcal{T}_G$ with crosses. The blue vertical lines are thus the nodal trajectories of $\mathcal{T}_0$ and the red skewed vertical lines those of $\mathcal{T}_G$.}
\label{fig_spacetime_slabs_scmm}
\end{figure}
%
%
In general we will use bar, i.e., $\bar{\cdot}$, to denote something related to space-time, such as domains and variables. In addition to the domains $\Omt{1}$ and $\Omt{2}$, we also consider the ``covered'' overlap domain $\Omt{O}$. To define it we will use the set of simplices $\mathcal{T}_{0, \bGn} := \{K \in \mathcal{T}_0 : \exists t \in I_n \text{ such that } K \cap \Gamma(t) \neq \emptyset \}$, i.e., all simplices in $\mathcal{T}_0$ that are cut by $\bGn$. We define the overlap domain $\Omt{O}$ for a time $t \in I_n$ by
\begin{equation}
\Omt{O} := \bigcup_{K \in \mathcal{T}_{0, \bGn}} K \cap \Omt{2}
\label{def:OmOt}
\end{equation}
As a discrete counterpart to the motion of the domain $G$, we prescribe a simple continuous motion for the overlapping mesh $\mathcal{T}_G$. By this we mean that the location of the overlapping mesh $\mathcal{T}_G$ is a function with respect to time that is \emph{continuous} on $[0, T]$ and \emph{linear} on each $I_n$. This means that the discrete velocity we prescribe for $\mathcal{T}_G$ is constant on each $I_n$. Henceforth, we let $\mu$ denote this discrete velocity. Letting $\mu_\text{cont}$ denote the velocity prescribed for $G$, we take the discrete velocity to be $\mu|_{I_n} = k_n^{-1} \int_{I_n} \mu_\text{cont}(t) \ud t$, for $n = 1, \dots, N$, i.e., the slabwise average. An illustration of the slabwise space-time domains $S_{i, n}$ defined by \eqref{def:Sin} is shown in Figure~\ref{fig_spacetime_slabs_scmm} (Left). Figure~\ref{fig_spacetime_slabs_scmm} (Right) shows a slabwise space-time discretization that has both straight and skewed space-time trajectories as a result of the simple continuous mesh motion. In a standard setting with only straight space-time trajectories, the time-derivative operator $\partial_t$ is naturally also a derivative operator in the direction of the trajectories. This is convenient and we would like have an analogous operator for our setting. We start by defining the domain-dependent velocity $\mu_i = \mu_i(t)$ by
\begin{equation}
\mu_i(t) := 
\left\{
\begin{split}
\mathbf{0} \quad i = 1 \\
\mu(t) \quad i = 2
\end{split}
\right.
\label{def:mui}
\end{equation}
We use this velocity to define the domain-dependent derivative operator $D_t = D_{t,i}$ by
\begin{equation}
D_{t, i}\{\cdot\} := \partial_{t}\{\cdot\} + \mu_i \cdot \nab\{\cdot\}
\label{def:Dt}
\end{equation}
The operator $D_t$ is a \emph{scaled} derivative operator in the direction of the space-time trajectories. To see this, consider the space-time vector $\bm_i = (\mu_i, 1)$ and the space-time gradient $\bnab = (\nab, \partial_t)$. The \emph{unscaled} derivative operator in the direction of the space-time trajectories is
\begin{equation}
D_{s, i} = \frac{\bm_i}{|\bm_i|} \cdot \bnab = \frac{1}{|\bm_i|} \big( \mu_i \cdot \nab + \partial_t \big) =  \frac{1}{|\bm_i|} D_{t, i}
\label{}
\end{equation}
We thus have $D_{t, i} = |\bm_i| D_{s, i}$. Let $\bar{\tau} = \bar{\tau}(t)$ denote a space-time trajectory that is uncut on the time interval $(t_a, t_b)$, and $v$ be a function of sufficient regularity. The intrinsic scaling of $D_t$ gives the convenient integral identity 
\begin{equation}
 \int_{\bar{\tau}(t_a)}^{\bar{\tau}(t_b)} D_s v \ud s = \int_{t_a}^{t_b} D_t v \ud t
\end{equation}
\begin{wrapfigure}{hr}{0.5\textwidth}
\centering
\def\svgwidth{0.5\textwidth}
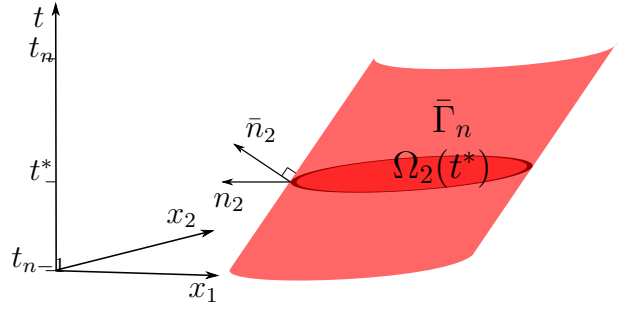
\caption{\emph{Space-time} normal vector $\bn_2$ to $\bGn$ (red) in relation to the \emph{spatial} normal vector $n_2$ to $\partial \Om{2}$.
\label{fig_spacetimevector}}
\vspace{-\baselineskip}
\end{wrapfigure}
Next we introduce some normal vectors. Let the spatial vector $n = n_i$ denote the outward pointing unit normal vector to $\partial \Omega_i$. Let the space-time vector $\bn = \bn_i = (\bn_i^x, \bn_i^t)$ denote the outward pointing unit normal vector to $\partial S_{i, n}$, where $\bn_i^x$ and $\bn_i^t$ denote the spatial and temporal component(s), respectively. On a purely \emph{spatial} subset, the space-time unit normal vector is purely \emph{temporal}, i.e., $\bn_i = (0, \pm 1)$, and vice versa, i.e., $\bn_i = (n_i, 0)$. The remaining case is a \emph{mixed} space-time subset and the only such set is $\bGn$. See Figure~\ref{fig_spacetimevector} for an illustration.
We define the space-time unit normal vector to $\bGn$ by
\begin{equation}
\bn_i |_{\bGn} = (\bn_i^x,  \bn_i^t)|_{\bGn} := \frac{1}{\sqrt{(n_i \cdot \mu)^2 + 1}} (n_i, -(n_i \cdot \mu))
\label{stvecg}
\end{equation}

\subsection{Finite element spaces} \label{subsec:fespaces}

\begin{figure}[h]
\centering
\def\svgwidth{0.5\textwidth}
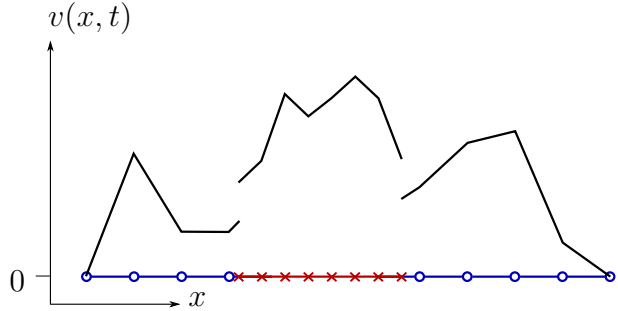
\caption{Example of $v(\cdot, t) \in V_h(t)$ for $d = 1$ and $p = 1$, where $\mathcal{T}_0$ is blue and $\mathcal{T}_G$ red.}
\label{figfefuncspace}
\end{figure}

We define the discrete spatial finite element spaces $V_{h,0}$ and $V_{h,G}$ as the spaces of continuous piecewise polynomials of degree $\le p$ on $\mathcal{T}_0$ and $\mathcal{T}_G$, respectively. We also let the functions in $V_{h,0}$ be zero on $\partial\Om{0}$. For $t \in [0, T]$, we use these two spaces to define the broken finite element space $V_h(t)$ by    
\begin{equation}
\begin{split}
V_h(t) := \{v : v|_{\Omt{1}} &= v_0|_{\Omt{1}} \text{ for some } v_0 \in V_{h,0} \text{ and } \\
v|_{\Omt{2}} &= v_G|_{\Omt{2}} \text{ for some } v_G \in V_{h,G} \}
\end{split} \label{fesvht}
\end{equation}
See Figure~\ref{figfefuncspace} for an illustration of a function $v \in V_h(t)$. 
\noindent For $n = 1, \dots, N$, we define the discrete space-time finite element spaces $V_{h,0}^n$ and $V_{h,G}^n$ as the spaces of functions that for a $t \in I_n$ lie in $V_{h,0}$ and $V_{h,G}$, respectively, and in time are polynomials of degree $\le q$ along the trajectories of $\mathcal{T}_0$ and $\mathcal{T}_G$ for $t \in I_n$, respectively. For $n = 1, \dots, N$, we use these two spaces to define the broken finite element space $V_h^n$ by:
\begin{equation}
\begin{split}
V_h^n := \{v : v|_{S_{1,n}} & = v_0^n|_{S_{1,n}} \text{ for some } v_0^n \in V_{h,0}^n \text{ and } \\ 
v|_{S_{2,n}} & = v_G^n|_{S_{2,n}} \text{ for some } v_G^n \in V_{h,G}^n \}
\end{split} \label{fesvhn}
\end{equation} 
We define the global space-time finite element space $V_h$ by: 
\begin{equation}
V_h := \{v : v|_{S_{0,n}} \in V_h^n, n = 1, \dots, N \} \label{fesvh}
\end{equation}

\subsection{Finite element formulation} \label{subsec:feform}

We may now formulate the space-time cut finite element formulation for the problem described in Section~\ref{sec:probform} as follows: Find $u_h \in V_h$ such that
\begin{equation}
B_h(u_h, v) = \int_{0}^T (f, v)_{\Om{0}} \ud t + (u_{0}, v_{0}^+)_{\Om{0}} \quad \forall v \in V_h
\label{feform}
\end{equation}
The non-symmetric bilinear form $B_h$ is defined by
\begin{equation}
\begin{split}
B_h(w, v) := & \;  \sum_{i=1}^2 \sum_{n=1}^N \int_{I_n} (\dot w, v)_{\Omt{i}} \ud t + \sum_{n=1}^N \int_{I_n} A_{h,t}(w, v) \ud t \\
& + \sum_{n=1}^{N-1}([w]_{n}, v_{n}^+)_{\Om{0}} + (w_{0}^+, v_{0}^+)_{\Om{0}} + \sum_{n=1}^N \int_{\bGn} -\bn^t [w]v_\sigma \ud\bs
\end{split}
\label{Bhdef}
\end{equation}
where $( \cdot , \cdot )_{\Omega}$ is the $L^2(\Omega)$-inner product, $[v]_n$ is the jump in $v$ at time $t_n$, i.e., $[v]_n = v_n^+ - v_n^-$, $v_n^\pm = \lim_{\varepsilon \to 0+} v(x, t_n \pm \varepsilon)$. The last term in $B_h$ mimics the standard dG-time-jump term, but over $\bGn$. Here, $\bn$ is the space-time normal vector to $\bGn$ defined by \eqref{stvecg}, $[v]$ is the jump in $v$ over $\bGn$, i.e., $[v]= v_1 - v_2$, $v_i = \lim_{\varepsilon \to 0+} v(\bs - \varepsilon \bn_i)$, $\bs = (s,t)$. If $\bn = \bn_1$, we take $\sigma = \frac{1}{2}(3 + \sgn(\bn^t))$ and if $\bn = \bn_2$, we take $\sigma = \frac{1}{2}(3 - \sgn(\bn^t))$, where sgn is the sign function. These choices make it so that $\sigma$ always picks the limit on the positive (in time) side of $\bGn$. The symmetric bilinear form $A_{h,t}$ is defined by
\begin{equation}
\begin{split}
A_{h,t}(w, v) := & \; \sum_{i=1}^2 (\nab w, \nab v)_{\Omt{i}} - \abm(\langle \partial_{\bn^x} w \rangle, [v])_{\Gamma(t)} - \abm(\langle \partial_{\bn^x} v \rangle, [w])_{\Gamma(t)} \\
& + \abm(\gamma h_K^{-1}[w],[v])_{\Gamma(t)} + ([\nab w],[\nab v])_{\Omt{O}}
\end{split}
\label{Ahdef}
\end{equation}
where $\abm = \sqrt{|\mu|^2 + 1}$, $\langle v \rangle$ is a convex-weighted average of $v$ on $\Gamma$, i.e., $\langle v \rangle = \omega_1v_1 + \omega_2v_2$, where $\omega_1, \omega_2 \in [0, 1]$ and $\omega_1 + \omega_2 = 1$, $\partial_{\bn^x} v = \bn^x \cdot \nabla v$, $\gamma \geq 0$ is a stabilization parameter, $h_K = h_K(x) = h_{K_0}$ for $x \in K_0$, where $h_{K_0}$ is the diameter of simplex $K_0 \in \mathcal{T}_0$, and $\Omt{O}$ is the overlap domain defined by \eqref{def:OmOt}. The reason for including the factor $\abm$ in the $\Gamma(t)$ terms is that when considering spacetime, these terms should be on $\bGn$. Since $\abm$ is the skewed temporal scaling, we have that
\begin{equation}
\int_{I_n} \abm(w, v)_{\Gamma(t)} \ud t = \int_{\bGn} wv  \ud \bs
\label{Gipint}
\end{equation}

\section{Analytic preliminaries}

\subsection{The bilinear form $A_{h,t}$}
The space of $A_{h, t}$ is $H^{3/2 + \varepsilon}(\cup_i \Omt{i})$ where $\varepsilon > 0$ may be arbitrarily small.
%
%
Let $\Gamma_{K}(t) := K \cap \Gamma(t)$. We define the following two mesh-dependent norms:
\begin{equation}
\| w \|_{1/2,h,\Gamma(t)}^2 := \sum_{K \in \mathcal{T}_{0,\Gamma(t)}} h_K^{-1} \| w \|_{\Gamma_K(t)}^2 \qquad \| w \|_{-1/2,h,\Gamma(t)}^2 := \sum_{K \in \mathcal{T}_{0,\Gamma(t)}} h_K \| w \|_{\Gamma_K(t)}^2
\label{def:HHnorm} 
\end{equation}
Note that
\begin{equation}
\| w \|_{\Gamma(t)}^2 \leq h \| w \|_{1/2,h,\Gamma(t)}^2 \qquad (w,v)_{\Gamma(t)} \leq \| w \|_{-1/2,h,\Gamma(t)}\| v \|_{1/2,h,\Gamma(t)}
\label{HHnormineqs}
\end{equation}
We define the time-dependent spatial energy norm $\norma{\cdot}$ by
\begin{equation}
\norma{w}^2 := \sum_{i = 1}^2 \| \nab w \|_{\Omt{i}}^2 + \abm\|\langle \partial_{\bn^x} w \rangle \|_{-1/2,h,\Gamma(t)}^2 + \abm\|[w] \|_{1/2,h,\Gamma(t)}^2 + \|[\nab w]\|_{\Omt{O}}^2
\label{def:anorm}   
\end{equation}   
Continuity of $A_{h,t}$ follows from using \eqref{HHnormineqs} in \eqref{Ahdef}. Next we consider the coercivity: 
 
\begin{lemma}[Discrete coercivity of $A_{h,t}$] \label{Ahtcoerlem} 
Let the bilinear form $A_{h,t}$ and the energy norm $\norma{\cdot}$ be defined by \eqref{Ahdef} and \eqref{def:anorm}, respectively. Then, for $t \in [0, T]$ and $\gamma$ sufficiently large,
\begin{equation}
A_{h,t}(v,v) \gtrsim \norma{v}^2 \quad \forall v \in V_h(t) \label{Ahtcoer}
\end{equation}
\begin{proof} Following the proof of the coercivity in \cite{Hansbo:2002aa}, we consider
\begin{equation}
\begin{split}
2\abm(\langle \partial_{\bn^x} v \rangle, [v])_{\Gamma(t)} \le & \; \frac{\abm}{\varepsilon}\|\langle \partial_{\bn^x} v \rangle \|_{-1/2,h,\Gamma(t)}^2 + \varepsilon \abm \| [v] \|_{1/2,h,\Gamma(t)}^2 \\
\le & \; \frac{2\abm}{\varepsilon} C_I \bigg( \sum_{i=1}^2 \| \nab v \|_{\Omt{i}}^2 + \| [\nab v]\|_{\Omt{O}}^2 \bigg) \\
& - \frac{\abm}{\varepsilon}\|\langle \partial_{\bn^x} v \rangle \|_{-1/2,h,\Gamma(t)}^2 + \varepsilon \abm \| [v] \|_{1/2,h,\Gamma(t)}^2
\end{split} \label{Ahtcoerlemvvmid}
\end{equation}
where we have used Lemma~\ref{invineqgamlem} and denoted its constant by $C_I$ . We use \eqref{Ahtcoerlemvvmid} in
\begin{equation}
\begin{split}
A_{h,t}(v,v) = & \; \sum_{i=1}^2 \|\nab v \|_{\Omt{i}}^2 - 2\abm(\langle \partial_{\bn^x} v \rangle, [v])_{\Gamma(t)} + \gamma \abm \| [v]\|_{1/2,h,\Gamma(t)}^2 + \|[\nab v]\|_{\Omt{O}}^2 \\
\ge & \; \bigg(1 - \frac{2 \abm C_I}{\varepsilon}\bigg) \sum_{i=1}^2\|\nab v \|_{\Omt{i}}^2 + \frac{\abm}{\varepsilon}\|\langle \partial_{\bn^x} v \rangle \|_{-1/2,h,\Gamma(t)}^2 \\
& + (\gamma - \varepsilon ) \abm\| [v] \|_{1/2,h,\Gamma(t)}^2 + \bigg(1 - \frac{2 \abm C_I}{\varepsilon} \bigg) \|[\nab v]\|_{\Omt{O}}^2 
\end{split} \label{Ahtcoerlemvvfin}
\end{equation}
By taking $\varepsilon > 2\abm C_I$, and $\gamma > \varepsilon$ we may obtain (\ref{Ahtcoer}) from (\ref{Ahtcoerlemvvfin}).
\end{proof}
\end{lemma}

\subsection{The bilinear form $B_h$}

The bilinear form $B_h$ can be expressed differently, as noted in the following lemma: 

\begin{lemma}[Alternative form of $B_h$] \label{lem:Bhpit} 
Let $\zeta = \frac{1}{2}(3 - \emph{sgn}(\bn^t))$. The bilinear form $B_h$, defined by \eqref{Bhdef}, can be written as
\begin{equation}
\begin{split}
B_h(w, v) = & \; \sum_{i=1}^2 \sum_{n=1}^N \int_{I_n} (w,- \dot v)_{\Omt{i}} \ud t + \sum_{n=1}^N \int_{I_n} A_{h,t}(w, v) \ud t \\
& + \sum_{n=1}^{N-1}(w_{n}^-,-[v]_{n})_{\Om{0}} + (w_{N}^-, v_{N}^-)_{\Om{0}} + \sum_{n=1}^N \int_{\bGn} \bn^t w_{\zeta}[v] \ud\bs
\end{split} \label{Bhpit}
\end{equation}
\begin{proof}
The proof is analogous to the standard case. The first term in \eqref{Bhdef} is integrated by parts in time via $\int_{S_{i,n}} (\nab, \partial_t) \cdot (\mathbf{0}, w v) \ud \bx$ and the result is combined with the last three terms in \eqref{Bhdef}. The combination of purely time-jump-related terms is exactly as in the standard case. For the $\bGn$-integral terms, we let $\zeta = \frac{1}{2}(3 - \sgn(\bn^t))$, if $\sigma = \frac{1}{2}(3 + \sgn(\bn^t))$ and $\bn = \bn_1$. This makes $\zeta, \sigma \in \{1, 2\}$ and $\zeta \neq \sigma$.
\end{proof}
\end{lemma}

\noindent An important result for the analysis is obtained by first taking the same function as both arguments of $B_h$. We present this result as a coercivity of $B_h$ with the following space-time energy norm:
\begin{equation}
\begin{split}
\normb{v}^2 := & \; \sum_{n=1}^N \int_{I_n} \norma{v}^2 \ud t \\
& + \sum_{n=1}^{N-1} \|[v]_n \|_{\Om{0}}^2 + \| v_N^- \|_{\Om{0}}^2 + \| v_0^+ \|_{\Om{0}}^2 + \sum_{n=1}^N \| |\bn^t|^{1/2} [v]\|_{\bGn}^2
\end{split}
\label{def:bnorm}
\end{equation}
\begin{lemma}[Discrete coercivity of $B_h$] \label{lem:Bhcoer}
Let the bilinear form $B_h$ and the energy norm $\normb{\cdot}$ be defined by \eqref{Bhdef} and \eqref{def:bnorm}, respectively. Then, for $\gamma$ sufficiently large, we have that 
\begin{equation}
B_h(v, v) \gtrsim \normb{v}^2 \quad \forall v \in V_h \label{Bhcoer}
\end{equation}
\begin{proof}
The proof is analogous to the standard case. First the same function $v$ is taken as both arguments of $B_h$. Then the first term in \eqref{Bhdef} is integrated in time via $\int_{S_{i,n}} (\nab, \partial_t) \cdot (\mathbf{0}, v^2) \ud \bx$ and the result is combined with the last three terms in \eqref{Bhdef}. The combination of purely time-jump-related terms is exactly as in the standard case. For the $\bGn$-integral terms, we note from the interdependence of $\sigma$ and $\bn$ that the combined integrand may be written as $\bn^t \sgn(\bn^t)[v]^2$. Also using Lemma~\ref{Ahtcoerlem} then shows the desired estimate.
\end{proof}
\end{lemma}

\noindent For the continued analysis, we define three space-time energy norms by
\begin{align}
\normx{v}^2 := & \; \sum_{i=1}^2 \sum_{n=1}^N \int_{I_n} {k_n} \| D_t v \|_{\Omt{i}}^2 \ud t + \normb{v}^2
\label{def:xnorm} \\
\normyp{v}^2 := & \; \sum_{n=1}^N \bigg(\int_{I_n} \frac{1}{k_n} \| v \|_{\Om{0}}^2 \ud t + \int_{I_n} \norma{v}^2 \ud t + \|v_{n-1}^+ \|_{\Om{0}}^2 \bigg) \label{def:ypnorm} \\
\normym{v}^2 := & \; \sum_{n=1}^N \bigg(\int_{I_n} \frac{1}{k_n} \| v \|_{\Om{0}}^2 \ud t + \int_{I_n} \norma{v}^2 \ud t + \|v_n^- \|_{\Om{0}}^2 \bigg)
\label{def:ymnorm}
\end{align}
The $X$-norm is the main norm, meaning that it is in this norm that we obtain stability and error estimates. The $Y$-norms are auxiliary norms. We use the $X$-norm and $Y$-norms to obtain continuity of $B_h$ which comes in two variants depending on the starting point, i.e., the standard form of $B_h$ \eqref{Bhdef} or the alternative \eqref{Bhpit}. 
\begin{lemma}[Continuity of $B_h$] \label{lem:Bhcont}
Let the bilinear form $B_h$ be defined by \eqref{Bhdef} and the norms $\normx{\cdot}$, $\normyp{\cdot}$, and $\normym{\cdot}$ by \eqref{def:xnorm}, \eqref{def:ypnorm}, and \eqref{def:ymnorm}, respectively. Then for any functions $w$ and $v$ with sufficient spatial and temporal regularity we have that
\begin{align}
B_h(w, v) & \lesssim \normx{w} \normyp{v} \label{Bhcontp} \\
B_h(w, v) & \lesssim \normym{w} \normx{v} \label{Bhcontm}
\end{align} 
\begin{proof}
The proofs of \eqref{Bhcontp} and \eqref{Bhcontm} are analogous so we only consider the latter since it gives the continuity result needed in the error analysis. The starting point is the alternative form of $B_h$ \eqref{Bhpit}. Applying the Cauchy--Schwarz inequality to all the terms (several times and different versions for some), \eqref{def:Dt} to split the first term followed by Corollary~\ref{cor:poincare_ineq_energy} for the $w$-factor in the resulting $\mu_i \cdot \nab$-part, the continuity of $A_{h,t}$ in the treatment of the second term, and Lemma~\ref{lem:invineqGam_Vterm} in the treatment of the fifth, we get product terms, where one factor may be estimated by $\normym{w}$ and the other by $\normx{v}$.
\end{proof}
\end{lemma}

\noindent Next, we present an estimate involving the bilinear form $B_h$ and the $X$-norm that may be viewed as a counterpart to such a coercivity. Due to the appearance of the estimate, we call it ``perturbed coercivity''. The estimate is a cornerstone of the energy analysis. It is fundamental to the stability analysis and also the starting point for deriving an inf-sup condition that in turn is essential for the error analysis. Key technical results used in the proof of the perturbed coercivity are Lemma~\ref{lem:invest_temp_anorm} and Lemma~\ref{lem:inv_ineq_temp_om}.

\begin{lemma}[Discrete perturbed coercivity of $B_h$] \label{lem:Bhpercoer}
Let the bilinear form $B_h$ and the norm $\normx{\cdot}$ be defined by \eqref{Bhdef} and \eqref{def:xnorm}, respectively. Then, for $q = 0, 1$, and $\gamma$ sufficiently large, there exists a constant $\delta > 0$ such that 
\begin{equation}
B_h(v, v + \delta k_n D_t v) \gtrsim \normx{v}^2 \quad \forall v \in V_h
\label{lemres:Bhpercoer}
\end{equation}

\begin{proof}
Using Lemma~\ref{lem:Bhcoer} with constant $\beta > 0$, the left-hand side of \eqref{lemres:Bhpercoer} is
\begin{equation}
B_h(v, v + \delta k_n D_t v) \geq \beta \normb{v}^2 +  B_h(v, \delta k_n D_t v)
\label{Bhpercoer_split} 
\end{equation}
The second term on the right-hand side is
\begin{equation}
\begin{split}
B_h(v, \delta k_n D_t v) = & \sum_{i=1}^2 \sum_{n=1}^N \int_{I_n} (\dot v, \delta k_n D_t v)_{\Omt{i}} + \sum_{n=1}^N \int_{I_n} A_{h,t}(v, \delta k_n D_t v) \ud t \\
& + \sum_{i=1}^2\sum_{n=1}^{N-1}([v]_{n}, (\delta k_n D_t v)_{n}^+)_{\Om{{i,n}}} + \sum_{i=1}^2(v_{0}^+, (\delta k_n D_t v)_{0}^+)_{\Om{{i,0}}} \\
& + \sum_{n=1}^N \int_{\bGn} -\bn^t [v](\delta k_n D_t v)_\sigma \ud\bs  
\end{split}
\label{Bhpercoer_Bvkdv0}
\end{equation}
The treatment of most of the terms involve the Cauchy--Schwarz inequality and for some also an $\varepsilon$-weighted Young's inequality. The first term in \eqref{Bhpercoer_Bvkdv0} is split using \eqref{def:Dt}, where the $D_t$-part is good, and we use standard estimates for the $\mu_i \cdot \nab$-part. For the second term in \eqref{Bhpercoer_Bvkdv0}, we use the continuity of $A_{h,t}$ followed by Lemma~\ref{lem:invest_temp_anorm}. The third and fourth term in \eqref{Bhpercoer_Bvkdv0} are estimated by Lemma~\ref{lem:inv_ineq_temp_om}. For the fifth and final term in \eqref{Bhpercoer_Bvkdv0}, we use Lemma~\ref{lem:invineqGam_Vterm} and Lemma~\ref{lem:invest_temp_anorm}. Collecting all the estimates and using the result in \eqref{Bhpercoer_split}, we may obtain 
\begin{equation}
\begin{split}
B_h(v, v + \delta k_n D_t v) \geq & \; \delta\bigg(1 - \frac{\delta}{\varepsilon} C \bigg) \sum_{i=1}^2 \sum_{n=1}^N \int_{I_n} k_n\| D_t v \|_{\Omt{i}}^2 \ud t \\
& + \bigg(\beta - \bigg(\varepsilon + \delta + \frac{\delta^2}{\varepsilon}\bigg)C\bigg) \normb{v}^2
\end{split}
\label{Bhpercoer_split_collect} 
\end{equation}
where $C > 0$ denote various constants. First taking $\varepsilon > 0$ sufficiently small and then taking $\delta > 0$ sufficiently small gives the desired estimate.
\end{proof}
\end{lemma}

\noindent Using Lemma~\ref{lem:Bhpercoer} and Lemma~\ref{lem:invest_temp_xnorm}, we may obtain the discrete inf-sup condition:

\begin{corollary}[A discrete inf-sup condition for $B_h$] \label{cor:Bhinfsup}
Let the bilinear form $B_h$ and the norm $\normx{\cdot}$ be defined by \eqref{Bhdef} and \eqref{def:xnorm}, respectively.
Then, for $q = 0, 1$, and $\gamma$ sufficiently large, we have that 
\begin{equation}
\normx{w} \lesssim \sup_{v \in V_h \setminus \{0\}} \frac{B_h(w, v)}{\normx{v}} \quad \forall w \in V_h
\label{Bhinfsup_res}
\end{equation}
\end{corollary}

\noindent To show Galerkin orthogonality, we need the following lemma on consistency:

\begin{lemma}[Consistency] \label{cnstylem}
The solution $u$ to problem \eqref{heateq} also solves \eqref{feform}.  

\begin{proof}
First insert $u$ in place of $u_h$ on the left-hand side of \eqref{feform} and use the regularity of $u$. Then integrate by parts in space via $\int_{S_{i,n}} (\nab, \partial_t) \cdot (\nab u v, 0) \ud \bx$ to get interior and boundary terms. The exterior boundary terms vanish because of the boundary conditions imposed on $v$ thus leaving the $\Gamma$-terms which are combined. Applying Lemma~\ref{jmplem} and the regularity of $u$ only leaves terms which from \eqref{heateq} equals the right-hand side of \eqref{feform}.
\end{proof}
\end{lemma}

\noindent From Lemma~\ref{cnstylem}, we may obtain the Galerkin orthogonality:

\begin{corollary}[Galerkin orthogonality]\label{cor:galort}
Let the bilinear form $B_h$ be defined by \eqref{Bhdef}, and let $u$ and $u_h$ be the solutions of \eqref{heateq} and \eqref{feform}, respectively. Then
\begin{equation}
B_h(u - u_h, v) = 0 \quad \forall v \in V_h \label{galort}
\end{equation}
\end{corollary} 


\section{Stability analysis}\label{seb:stabilityanalysis}

In this section we present and prove a stability estimate for the solution $u_h$ to \eqref{feform}. The key component in the proof is Lemma~\ref{lem:Bhpercoer}, i.e., the perturbed coercivity of $B_h$ on $V_h$. 

\begin{lemma}[A stability estimate in $\normx{\cdot}$] \label{stablem1}
Let $u_h$ be the solution of \eqref{feform}. Let $u_0$ and $f$ be the initial data and source in \eqref{heateq}, respectively. Then, for $q = 0, 1$, and $\gamma$ sufficiently large, we have that
\begin{equation}
\normx{u_h} \lesssim \| u_{0} \|_{\Om{0}} + \| f \|_{L^2((0,T]; L^2(\Om{0}))}
\label{stablemres}
\end{equation}
\begin{proof}
By taking $v = u_h \in V_h$ in Lemma \ref{lem:Bhpercoer} and $v = u_h + \delta k_n D_t u_h \in V_h$ in \eqref{feform}, we have
\begin{equation}
\begin{split}
\normx{u_h}^2 \lesssim & \; B_h(u_h, u_h + \delta k_n D_t u_h) \\
= & \; (u_{0}, u_{h, 0}^+)_{\Om{0}} + (u_{0}, \delta k_1 (D_t u_h)_0^+)_{\Om{0}} \\
& + \sum_{n=1}^N \int_{I_n} (f, u_h)_{\Om{0}} \ud t + \sum_{n=1}^N \int_{I_n} (f, \delta k_n D_t u_h)_{\Om{0}} \ud t
\end{split}
\label{stablem0}
\end{equation}
Applying the Cauchy--Schwarz inequality to all the terms (several times and different versions for some), Lemma~\ref{lem:inv_ineq_temp_om} in the treatment of the second term, and Corollary~\ref{cor:poincare_ineq_energy} in the treatment of the third, we get product terms, where one factor is $\| u_{0} \|_{\Om{0}}$ or $\| f \|_{L^2((0,T]; L^2(\Om{0}))}$ and the other may be estimated by $\normx{u_h}$. Dividing both sides by $\normx{u_h}$ thus gives \eqref{stablemres}. 

\end{proof} 

\end{lemma}

\section{A priori error analysis}\label{sec:errorestimate}

\begin{theorem}[An optimal order a priori error estimate in $\normx{\cdot}$]
\label{thm:errorest}
Let $\normx{\cdot}$ be defined by $(\ref{def:xnorm})$, let $u$ be the solution of \eqref{heateq} and let $u_h$ be the finite element solution defined by \eqref{feform}. Then, for $q = 0, 1$, and $\gamma$ sufficiently large, we have that
\begin{equation}
\normx{u - u_h}^2 \lesssim k^{2q+1}F_k^2(u) + h^{2p} \bigg(F_h^2(u) + E_{h, 1}^2(u) \bigg)
\label{errorest_res}
\end{equation}
where $F_k$, $F_h$, and $E_{h,1}$ are defined by \eqref{interpest_res_Fk}, \eqref{interpest_res_Fh}, and \eqref{interpest_glob_st_Eh}, respectively.
\begin{proof}
We use the interpolant $\interph u \in V_h$, where $\interph$ is the space-time interpolation operator defined by (\ref{def:interph}), to split the error $e = u - u_h$ into $\rho = u - \interph u$ and $\theta = \interph u - u_h$. Thus
\begin{equation}
\normx{e} \leq \normx{\rho} + \normx{\theta}
\label{errorsplit_norm}
\end{equation}
where we focus on the $\theta$-part first. From Corollary~\ref{cor:galort}, i.e., Galerkin orthogonality, we have for any $v \in V_h$ that
\begin{equation}
B_h(\theta, v) = -B_h(\rho, v)
\label{errorrep_theta_galort}
\end{equation}
We note that $\theta \in V_h$ and use Corollary~\ref{cor:Bhinfsup}, i.e., a discrete inf-sup condition for $B_h$, the Galerkin orthogonality result (\ref{errorrep_theta_galort}), and Lemma~\ref{lem:Bhcont}, i.e., continuity of $B_h$, to estimate the $\theta$-part by
\begin{equation}
\begin{split}
\normx{\theta} \lesssim \sup_{v \in V_h \setminus \{0\}} \frac{B_h(\theta, v)}{\normx{v}} = \sup_{v \in V_h \setminus \{0\}} \frac{-B_h(\rho, v)}{\normx{v}} \lesssim \sup_{v \in V_h \setminus \{0\}} \frac{\normym{\rho} \normx{v}}{\normx{v}} = \normym{\rho}
\end{split}
\label{errorrep_theta_norm}
\end{equation}
Using (\ref{errorrep_theta_norm}) in (\ref{errorsplit_norm}), we estimate the approximation error by
\begin{equation}
\begin{split}
\normx{e}^2 \lesssim & \; \normx{\rho}^2 + \normym{\rho}^2 \\
\lesssim & \; \sum_{i,n} \bigg( k_n \int_{I_n} \| D_t \rho \|_{\Omt{i}}^2 \ud t + \frac{1}{k_n} \int_{I_n} \| \rho \|_{\Omt{i}}^2 \ud t \bigg) + \normb{\rho}^2 + \sum_{n=1}^N \|\rho_n^- \|_{\Om{0}}^2
\end{split}
\label{error_xnorm_est0}
\end{equation}
By applying various interpolation error estimates: Lemma~\ref{lem:interpest_glob_st_om0T} and using \eqref{quasiuniformity_st} for the first term, Lemma~\ref{lem:interpest} for the second, and Corollary~\ref{cor:interpest_glob_s} for the third, we get results that may be estimated by the right-hand side of \eqref{errorest_res}. 

\end{proof}

\end{theorem}


\section{Numerical results} \label{secnumres}

Here we present numerical results for a problem in one spatial dimension on the unit interval with exact solution $u(x, t) = \sin^2(\pi x)\euler^{-t/2}$. We compute $u_h$ for $p=1$ and $q = 0,1$. For dG(1) in time, \emph{some} of the left-hand side integrals involving time have been approximated locally by quadrature. For integrals over cut space-time prisms, \emph{composite three}-point Lobatto quadrature has been used in time. For integrals over intraprismatic segments of the space-time boundary $\bGn$, \emph{three}-point Lobatto quadrature has been used. Both of these choices of quadrature result in a quadrature error $= O(k^4)$. The right-hand side integrals have been approximated locally by quadrature over the space-time prisms: first quadrature in time, then quadrature in space. In space, the trapezoidal rule has been used, thus resulting in a quadrature error $= O(h^2)$. For dG(0) in time, the midpoint rule has been used, thus resulting in a quadrature error $= O(k^2)$. For dG(1) in time, \emph{three}-point Lobatto quadrature has been used, thus resulting in a quadrature error $= O(k^4)$. For simplicity, the velocity $\mu$ of the overlapping mesh is set to be constant at the value $\mu(t_n)$ on every subinterval $I_n = (t_{n-1}, t_n]$. The stabilization parameter $\gamma = 10$.

\begin{figure}[h]
\centering
\includegraphics[width=0.4\textwidth]{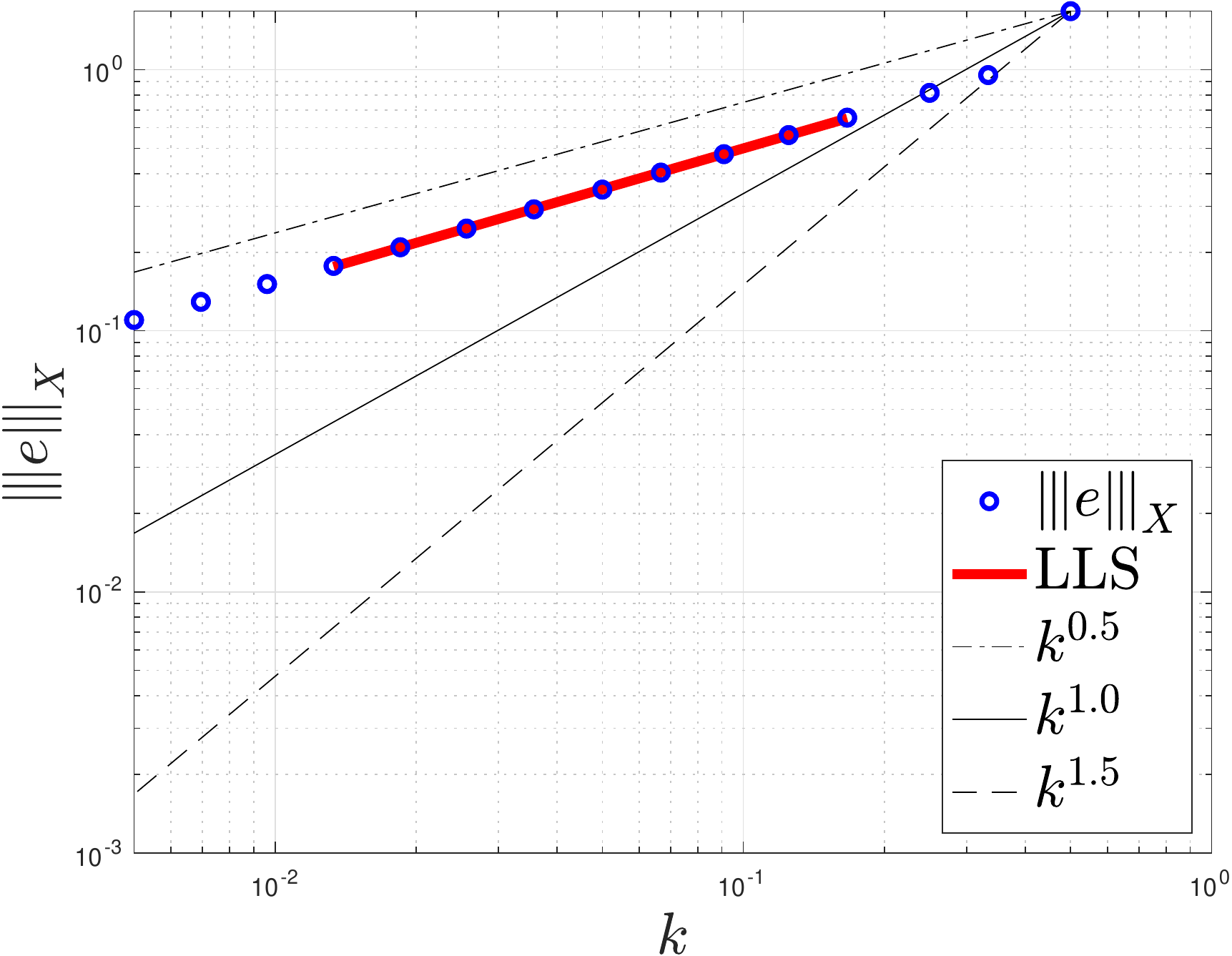}
\includegraphics[width=0.4\textwidth]{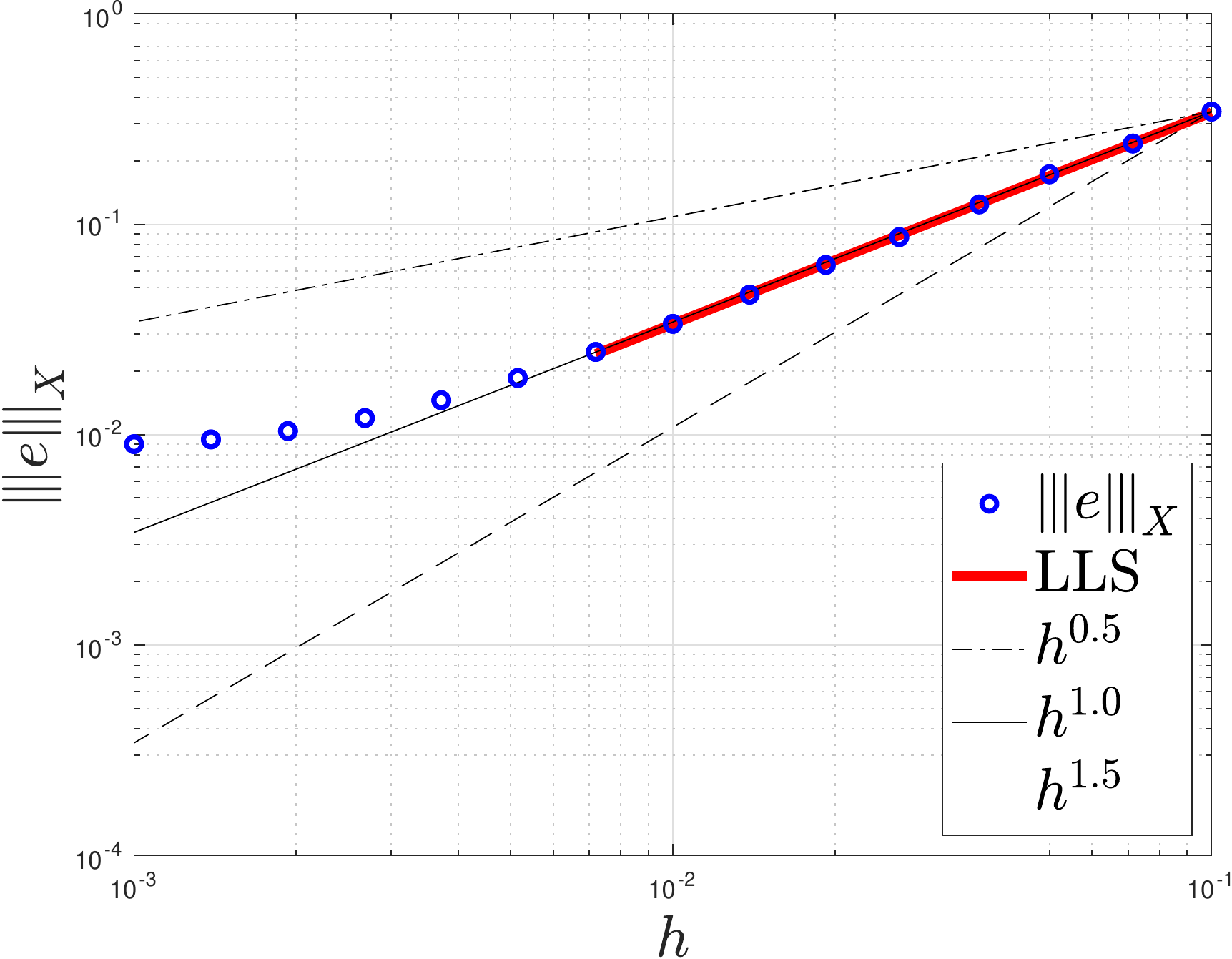}
\caption{Error convergence for dG(0) with $\mu = 0.6$.
\label{fig_dG0_ECC_mu0p6}}
\end{figure}

\begin{figure}[h]
\centering
\includegraphics[width=0.4\textwidth]{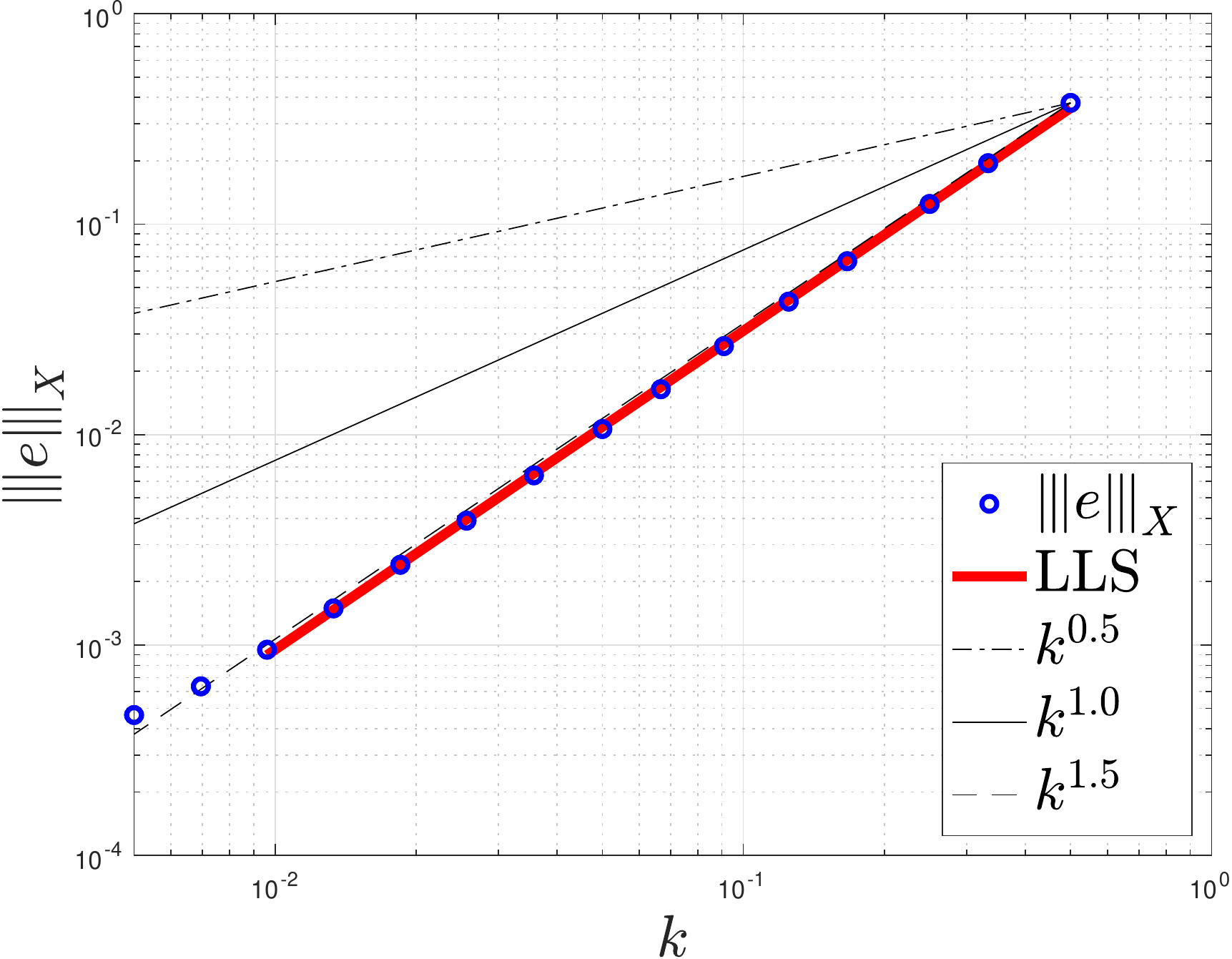}
\includegraphics[width=0.4\textwidth]{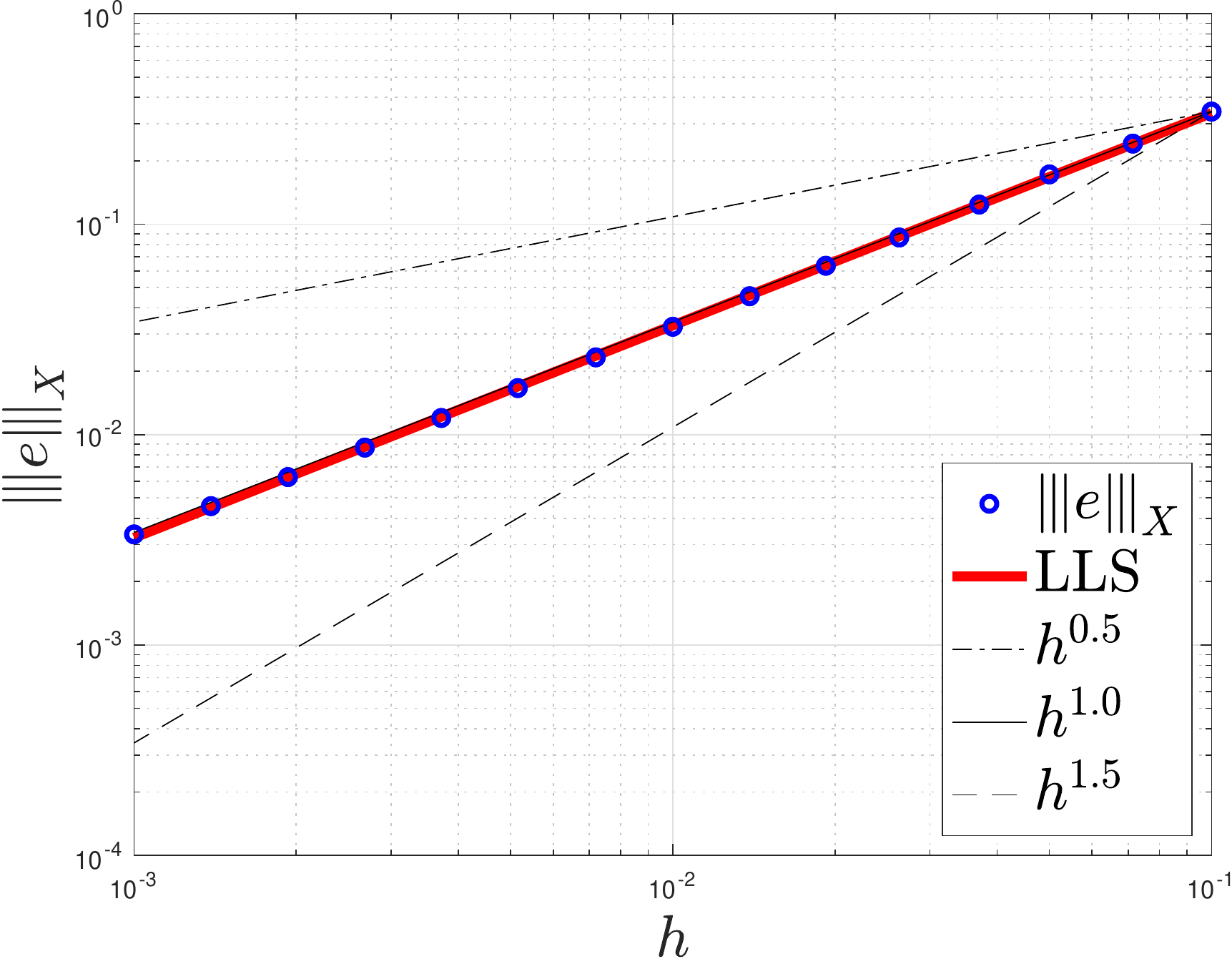}
\caption{Error convergence for dG(1) with $\mu = 0.6$.
\label{fig_dG1_ECC_mu0p6}}
\end{figure}

For the error convergence study, both $\mathcal{T}_0$ and $\mathcal{T}_G$ are uniform meshes, with mesh sizes $h_0$ and $h_G$, respectively. The temporal discretization is also uniform with time step $k$ for each instance. The final time is set to $T = 1$, the length of $\mathcal{T}_G$ is 0.25, and the initial position of $\mathcal{T}_G$ is the spatial interval [0.125, 0.125 + 0.25]. The $error$ is $\normx{e} = \normx{u - u_h}$. All time, space, and space-time integrals involving $u$ in the $X$-norm have been approximated locally by \emph{three}-point Gauss-Legendre quadrature: first quadrature in time, then quadrature in space where applicable. This results in a quadrature error $= O((k^6 + h^6)^{1/2})$. In the $k$-convergence study, the mesh sizes have been fixed at $h = 10^{-3}$ and $h = 10^{-4}$ for dG(0) and dG(1), respectively. Analogously, in the $h$-convergence study, the time step has been fixed at $k = 10^{-4}$ and $k = 10^{-2}$ for dG(0) and dG(1), respectively. Figure~\ref{fig_dG0_ECC_mu0p6} and \ref{fig_dG1_ECC_mu0p6} display error convergence plots for dG(0) and dG(1) in time with $\mu=0.6$. The left plots show the $error$ versus $k$ and the right plots versus $h = h_0 \geq h_G$. Besides the computed $error$, each plot contains a line segment that has been computed with the linear least squares method to fit the error data. This line segment is referred to as the LLS of the $error$. Reference slopes are also included. In Table~\ref{tabnumord} we summarize the slope of the LLS of the $error$ for different values of $\mu$. 
\begin{table}[h]
\centering
\begin{tabular}{ l | c | c | c | c | } 
  & \multicolumn{2}{|c|}{dG(0) in time} & \multicolumn{2}{|c|}{dG(1) in time} \\  
  \hline                    
  $\mu$ & versus $k$ (points) & versus $h$ (points) & versus $k$ (points) & versus $h$ (points) \\
  \hline
  0 & 0.5058 (4--12) & 1.0190 (1--9) & 1.4506 (1--5) & 1.0137 (1--15) \\
  \hline
  0.1 & 0.5026 (4--12) & 1.0205 (1--9) & 1.4893 (1--7) & 1.0155 (1--15) \\ 
  \hline
  0.2 & 0.4985 (4--12) & 1.0221 (1--9) & 1.4947 (1--8) & 1.0162 (1--15) \\ 
  \hline 
  0.4 & 0.5163 (4--12) & 1.0179 (1--9) & 1.5031 (1--11) & 1.0147 (1--15) \\ 
  \hline 
  0.6 & 0.5179 (4--12) & 1.0047 (1--9) & 1.5151 (1--13) & 1.0091 (1--15) \\
  \hline 
\end{tabular}
\caption{The slope of the LLS of the $error$ versus $k$ and $h$ for different values of $\mu$.
\label{tabnumord}}
\end{table}

\begin{figure}[h]
\centering
\includegraphics[width=0.32\textwidth]{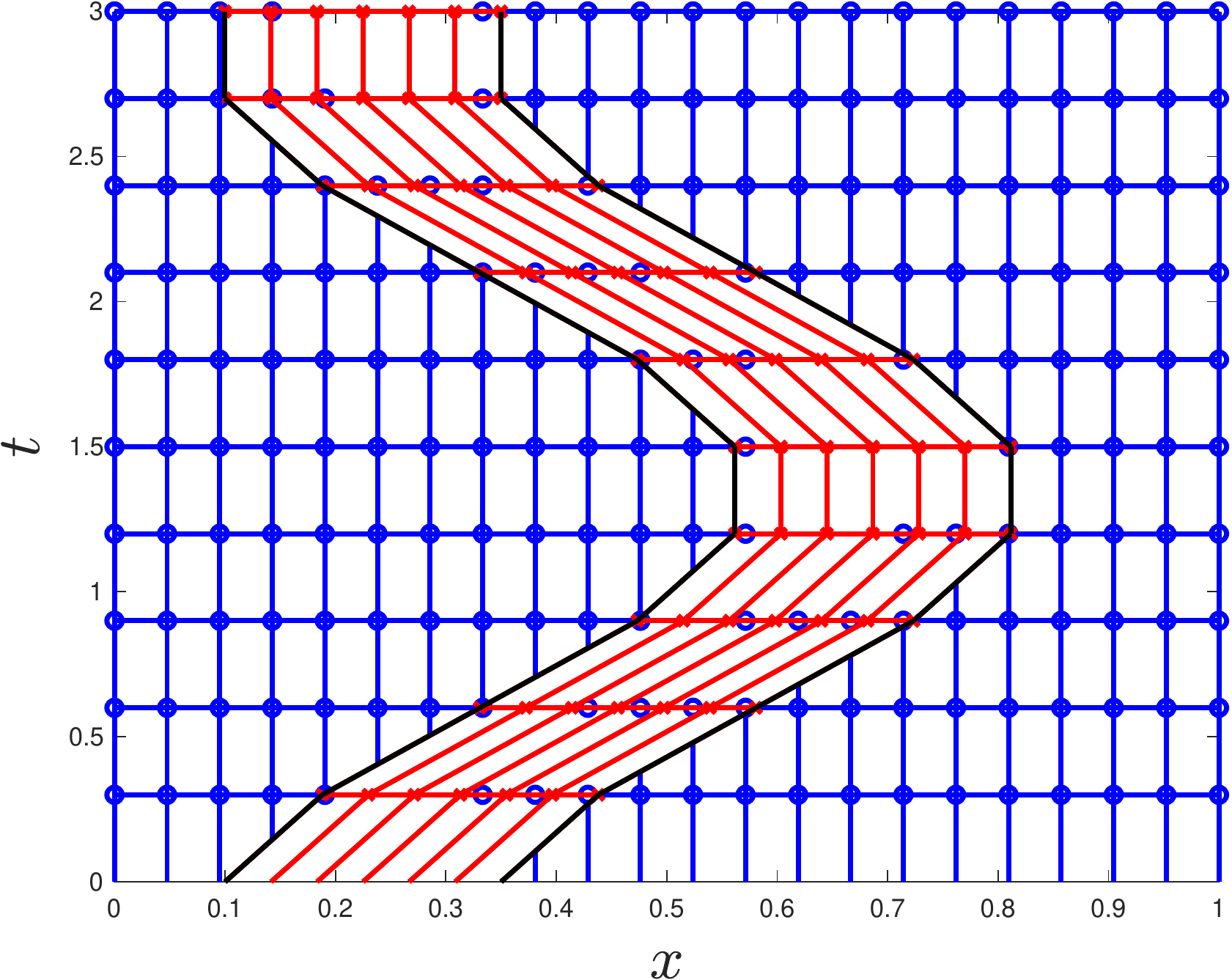}
\includegraphics[width=0.32\textwidth]{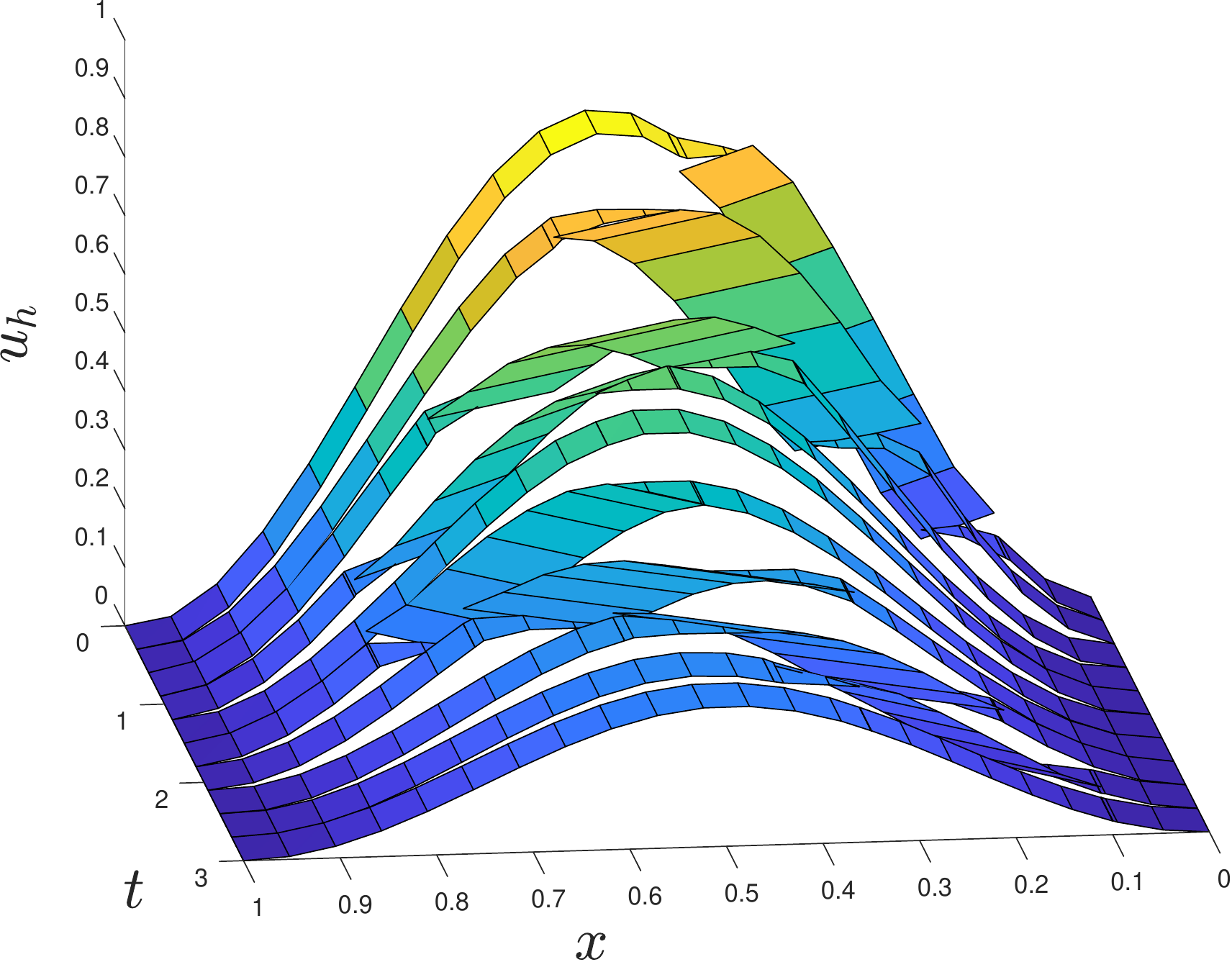}
\includegraphics[width=0.32\textwidth]{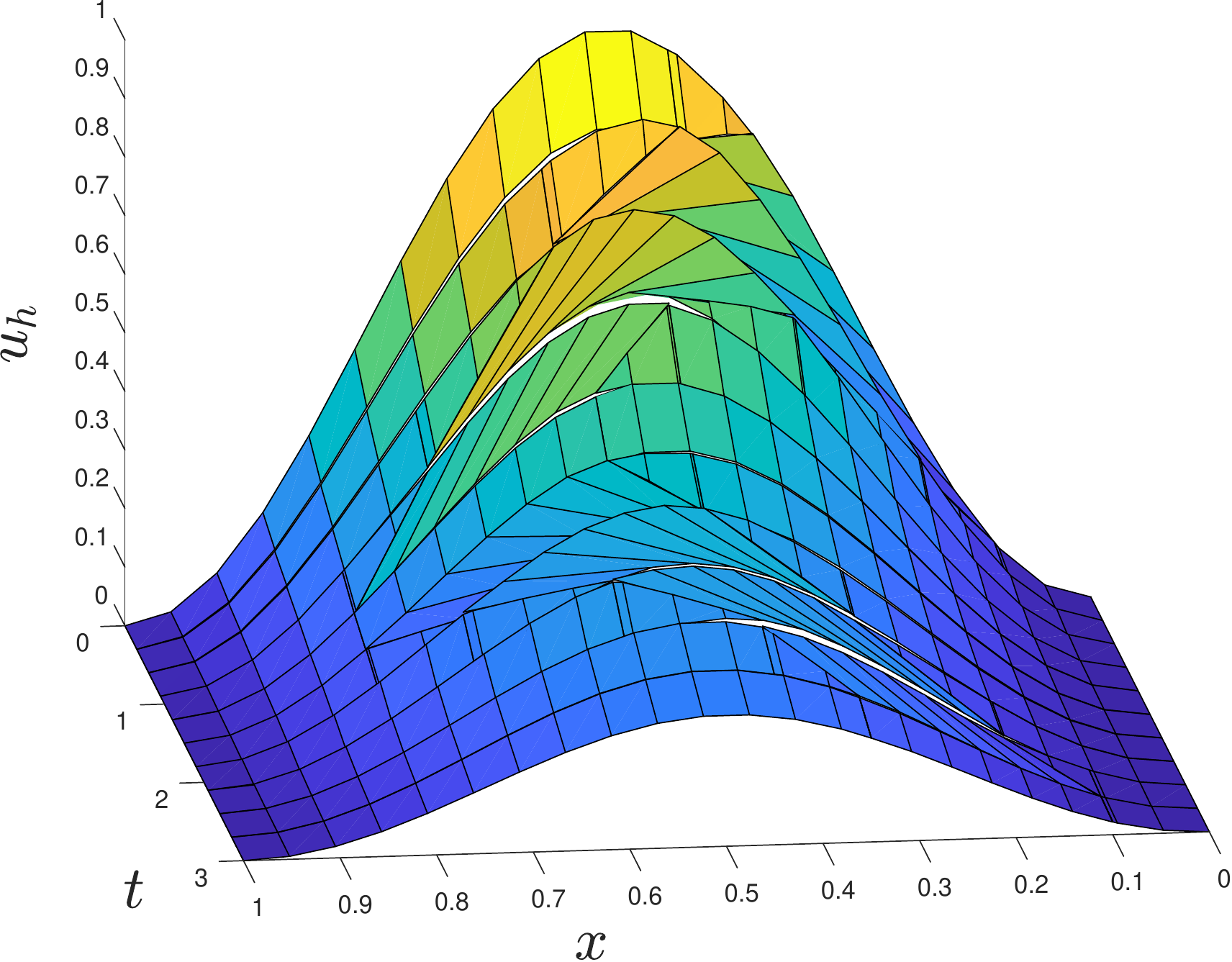}
\caption{Space-time discretization (left) with resulting dG(0)cG(1)-solution (middle) and dG(1)cG(1)-solution (right). \label{fig_numsol}}
\end{figure}

The numerical solutions presented in Figure~\ref{fig_numsol} have been computed for an equidistant space-time discretization: 22 nodes for $\mathcal{T}_0$, 7 nodes for $\mathcal{T}_G$ for all times, and 10 time steps on the interval $(0, 3]$. The length of $\mathcal{T}_G$ has again been 0.25 and the velocity $\mu$ has for simplicity been slabwise constant at $\mu|_{I_n} = \frac{1}{2}\sin(\frac{2 \pi t_n}{3})$.

\section{Conclusions}

We have presented a cut finite element method for a parabolic model problem on an overlapping mesh situation: one stationary background mesh and one \emph{continuously} moving overlapping mesh. We have applied what we believe to be  a relatively new analysis framework for finite element methods for parabolic problems. This new analysis framework may arguably be considered more robust and natural than standard ones, since it is the only one that we have been able to successfully apply to our overlapping mesh situation. The analysis is of an energy type and the main results are a basic stability estimate and an optimal order a priori error estimate. We have also presented numerical results for a parabolic problem in one spatial dimension that verify the analytic error convergence orders.


\appendix
\section{Analytic tools}

\begin{lemma}[A jump identity] \label{jmplem}
Let $\omega_+, \omega_- \in \mathbb{R}$ and $\omega_+ + \omega_- = 1$, let $[A] := A_+ - A_-$, and $\langle A \rangle := \omega_+A_+  + \omega_-A_-$. We then have
\begin{equation}
[AB] =  [A]\langle B \rangle + \langle A \rangle [B] + (\omega_- - \omega_+)[A][B]. \label{jmplemres}
\end{equation}
\begin{proof}
Using the definitions and evaluating both sides shows the identity.
\end{proof}
\end{lemma}

\subsection{Spatial estimates}

%
\begin{lemma}[A Poincar\'e inequality for $H_0^1(\cup_i \Omt{i})$]\label{lem:poincare_ineq_disc}
For $t \in [0, T]$ we have that
\begin{equation}
\| v \|_\Om{0} \lesssim \| \nab v \|_{\Omt{1}\cup\Omt{2}} + \|[v]\|_{\Gamma(t)} \quad \forall v \in H_0^1(\cup_i \Omt{i})
\label{poincare_ineq_res}
\end{equation}
\begin{proof}
To lighten the notation we omit the time dependence, which has no importance here anyways. For $v \in H_0^1(\Om{1} \cup \Om{2})$, we consider the dual problem: Find $\phi \in H^2(\Om{0}) \cap H_0^1(\Om{0})$ such that $- \lap{\phi} = v$ in $\Om{0}$. By using the dual problem, partial integration, that $v|_{\partial\Om{0}} = 0$, Lemma~\ref{jmplem}, and the regularity of $\phi$ ($[\partial_n\phi]|_\Gamma = 0$ in $L^2(\Gamma)$), we have
\begin{equation}
\| v \|_\Om{0}^2 = \sum_{i=1}^2(-\lap \phi, v)_{\Om{i}} = \sum_{i=1}^2 (\nab\phi, \nab v)_\Om{i} - (\langle \partial_n\phi \rangle, [v])_\Gamma
\label{poincare_ineq_pi}
\end{equation}
Using a standard trace inequality for $\nab\phi|_\Om{i} \in H^{1}(\Om{i})$, elliptic regularity on $H^2(\Om{0}) \cap H_0^1(\Om{0})$ for $\phi$, and the dual problem, the first argument to the last inner product may be estimated by
\begin{equation}
\|\langle \partial_n\phi \rangle \|_\Gamma \leq \sum_{i=1}^2\| \nab \phi_i \|_{\partial\Om{i}} \lesssim \sum_{i=1}^2 \| \nab \phi \|_{1,\Om{i}} \lesssim \| \phi \|_{2,\Om{0}} \lesssim \| \lap \phi \|_{\Om{0}} = \| v \|_{\Om{0}}
\label{poincare_ineq_avgest}
\end{equation} 
We note that this also gives an estimate for the first argument to the penultimate inner product. Thus using \eqref{poincare_ineq_avgest} in \eqref{poincare_ineq_pi} followed by cancellation of a factor $ \| v \|_\Om{0}$ on both sides gives \eqref{poincare_ineq_res}.
\end{proof}
\end{lemma}
\noindent By squaring both sides of (\ref{poincare_ineq_res}), using Young's inequality, and \eqref{HHnormineqs}, we may estimate the resulting right-hand side by $\norma{\cdot}$:
\begin{corollary}[An energy Poincar\'e inequality for $H^{3/2 + \varepsilon}(\cup_i \Omt{i}) \cap H_0^1(\cup_i \Omt{i})$]\label{cor:poincare_ineq_energy}
Let the time-dependent spatial energy norm $\norma{\cdot}$ be defined by \eqref{def:anorm}. Then, for $t \in [0, T]$, we have that
\begin{equation}
\| v \|_\Om{0} \lesssim \norma{v} \quad \forall v \in H^{3/2 + \varepsilon}(\cup_i \Omt{i}) \cap H_0^1(\cup_i \Omt{i})
\label{poincare_ineq_energy_res}
\end{equation}
\end{corollary}

\begin{lemma}[A spatial continuity result for $\Gamma(t)$] \label{lem:invineqGam_Vterm}
Let the space-time vector $\bn = (\bn^x, \bn^t)$ and the time-dependent spatial energy norm $\norma{\cdot}$ be defined by \eqref{stvecg} and \eqref{def:anorm}, respectively. Let $\sigma$ change arbitrarily along $\Gamma(t)$ between the values $1$ and $2$ and let $|\mu|_{[0, T]} = \max_{t \in [0, T]}\{ |\mu| \}$. Then, for $t \in [0, T]$, we have that
\begin{equation}
\begin{split}
(\bn^t [w], v_{\sigma})_{\Gamma(t)} \lesssim |\mu|_{[0, T]} h^{1/2} \norma{w} \norma{v} \\
\forall w, v \in H^{3/2 + \varepsilon}(\cup_i \Omt{i}) \cap H_0^1(\cup_i \Omt{i})
\end{split}
\label{lemres:invineqGam_Vterm}
\end{equation} 
\begin{proof}
To lighten the notation we omit the time dependence, which has no importance here anyways. Using $|\bn^t| \leq |\mu|$, which follows from (\ref{stvecg}), the left-hand side of (\ref{lemres:invineqGam_Vterm}) is
\begin{equation}
(\bn^t [w], v_{\sigma})_{\Gamma} \leq |\mu|_{[0, T]} \| [w]\|_{\Gamma} \|v_\sigma \|_{\Gamma}
\label{invineqGam_Vterm_0}
\end{equation}
Using \eqref{HHnormineqs} and \eqref{def:anorm}, the $w$-factor may be estimated by $ h^{1/2} \norma{w}$. Applying the standard trace inequality for $H^1(\Om{i})$, Corollary~\ref{cor:poincare_ineq_energy}, and \eqref{def:anorm}, the $v$-factor may be estimated by $\norma{v}$. This shows \eqref{lemres:invineqGam_Vterm}.
\end{proof}
\end{lemma}

\begin{lemma}[A scaled trace inequality for domain-partitioning manifolds of codimension 1] \label{lem:scatraineqdomcut}
For $d = 1, 2$, or $3$, let $\Omega \subset \mathbb{R}^d$ be a bounded domain with diameter $L$, i.e., $L = \text{diam}(\Omega) = \sup_{x, y \in \Omega} |x - y|$. Let $\Gamma \subset \Omega$ be a continuous manifold of codimension 1 that partitions $\Omega$ into $N$ subdomains. Then
\begin{equation}
\| v \|_{\Gamma}^2 \lesssim L^{-1}\| v \|_{\Omega}^2 + L \| \nab v \|_{\Omega}^2 \quad \forall v \in H^1(\Omega)
\label{lemres:scatraineqdomcut}
\end{equation}
\begin{proof}
If (\ref{lemres:scatraineqdomcut}) holds for the case $N = 2$, then that result may be applied repeatedly to show (\ref{lemres:scatraineqdomcut}) for $N > 2$. We thus assume that $\Gamma$ partitions $\Omega$ into two subdomains denoted $\Om{1}$ and $\Om{2}$ with diameters $L_1$ and $L_2$, respectively. From the regularity assumptions on $v$, we have for $i = 1, 2$, that $v \in H^1(\Om{i})$ and thus 
\begin{equation}
\| v \|_{\Gamma}^2 \leq \| v \|_{\partial \Om{i}}^2 \lesssim L_i^{-1}\| v \|_{\Om{i}}^2 + L_i \| \nab v \|_{\Om{i}}^2
\label{scatraineqdomcut_standard}
\end{equation}
where we have used a standard scaled trace inequality. Using the triangle type inequality $L \leq L_1 + L_2$ and \eqref{scatraineqdomcut_standard}, the left-hand side of \eqref{lemres:scatraineqdomcut} is
\begin{equation}
\| v \|_{\Gamma}^2 \leq \sum_{i=1}^2 \frac{L_i}{L}\| v \|_{\Gamma}^2 \lesssim \sum_{i=1}^2 \bigg(L^{-1}\| v \|_{\Om{i}}^2 + L \| \nab v \|_{\Om{i}}^2 \bigg) \lesssim L^{-1}\| v \|_{\Omega}^2 + L\| \nab v \|_{\Omega}^2
\label{scatraineqdomcut_0}
\end{equation}
which shows \eqref{lemres:scatraineqdomcut}.
\end{proof}
\end{lemma}
\noindent Let $\Gamma_K = \Gamma_K(t) = K \cap \Gamma(t)$. For $t \in [0, T]$, $j \in \{0, G\}$, a simplex $K \in \mathcal{T}_{j,\Gamma(t)} = \{K \in \mathcal{T}_j : K \cap \Gamma(t) \neq \emptyset\}$, and $v \in H^1(K)$, we have from Lemma~\ref{lem:scatraineqdomcut} that
\begin{equation}
\| v \|_{\Gamma_K}^2 \lesssim h_K^{-1}\| v \|_{K}^2 + h_K \| \nab v \|_{K}^2
\label{scatraineqGamK_warmup}
\end{equation}
where $h_K$ is the diameter of $K$. For $v \in \mathcal{P}(K)$, we have the standard inverse estimate
\begin{equation}
\| D_x^{k} v \|^2_K \lesssim h_K^{-2} \| D_x^{k-1} v \|^2_K \quad \text{ for } k \geq 1
\label{investpolK_standard}
\end{equation}
Using \eqref{investpolK_standard} in \eqref{scatraineqGamK_warmup}, we get the following corollary:
\begin{corollary}[A discrete spatial local inverse inequality for $\Gamma_K(t)$] \label{cor:scatraineqGamK}
For $t \in [0, T]$, $j \in \{0, G\}$, $K \in \mathcal{T}_{j,\Gamma(t)}$ with diameter $h_K$, let $\Gamma_K(t) = K \cap \Gamma(t)$. Then, for $k \geq 0$, we have that
\begin{equation}
\| D_x^k v \|_{\Gamma_K(t)}^2 \lesssim h_K^{-1} \| D_x^k v  \|_{K}^2 \quad \forall v \in V_h(t)
\label{corres:scatraineqGamK}
\end{equation}
\end{corollary}

\begin{lemma}[A discrete spatial inverse inequality for $\Gamma(t)$] \label{invineqgamlem}
Let the mesh-dependent norm $\| \cdot \|_{-1/2,h,\Gamma(t)}$ be defined by \eqref{def:HHnorm}. Then, for $t \in [0, T]$, we have that 
\begin{equation}
\| \langle \partial_{\bn^x} v \rangle \|_{-1/2,h,\Gamma(t)}^2 \lesssim \sum_{i=1}^2 \|\nab v\|_{\Omt{i}}^2 + \| [\nab v] \|_{\Omt{O}}^2 \quad \forall v \in V_h(t)
\label{invineqgamlemres}
\end{equation} 
\begin{proof}
To lighten the notation we omit the time dependence, which has no importance here anyways. We follow the proof of the corresponding inequality in \cite{Hansbo:2002aa} with some modifications. We use index $j \in \{0, G \}$, such that, if $j = 0$, then $i = 1$ and if $j = G$, then $i = 2$, and let $\Gamma_{K_j} = K_j \cap \Gamma$ and $\mathcal{T}_{j,\Gamma} = \{K_j \in \mathcal{T}_j : K_j \cap \Gamma \neq \emptyset\}$. Note that for $i=1, 2$,
\begin{equation}
\sum_{K_0 \in \mathcal{T}_{0,\Gamma}} h_{K_0} \| v_i \|_{\Gamma_{K_0}}^2 \lesssim \sum_{K_G \in \mathcal{T}_{G,\Gamma}} h_{K_G} \| v_i \|_{\Gamma_{K_G}}^2
\label{gamk0ineqgamkg}
\end{equation}
which follows from $\cup_{K_0 \in \mathcal{T}_{0,\Gamma}} \Gamma_{K_0} = \Gamma = \cup_{K_G \in \mathcal{T}_{G,\Gamma}} \Gamma_{K_G}$ and the inter-quasi-uniformity of the meshes. Since $\partial_{\bn^x} v = \bn^x \cdot \nab v$ and $|\omega_i| |\bn^x| \leq 1$, we have $\| \omega_i (\partial_{\bn^x} v)_i \|_{\Gamma_{K_j}}^2 \leq \|(\nab v)_i \|_{\Gamma_{K_j}}^2$. Using this after \eqref{gamk0ineqgamkg}, and followed by Corollary~\ref{cor:scatraineqGamK}, the left-hand side of \eqref{invineqgamlemres} is
\begin{equation}
\begin{split}
\| \langle \partial_{\bn^x} v \rangle \|_{-1/2,h,\Gamma}^2 & \lesssim \sum_{i=1}^2 \sum_{K_j \in \mathcal{T}_{j, \Gamma}} h_{K_j} \|(\nab v)_i \|_{\Gamma_{K_j}}^2 \lesssim \sum_{i=1}^2 \sum_{K_j \in \mathcal{T}_{j, \Gamma}} \|(\nab v)_i\|_{K_j}^2 \\
& = \sum_{K_0 \in \mathcal{T}_{0,\Gamma}} \bigg( \|\nab v\|_{K_0 \cap \Om{1}}^2 +  \|(\nab v)_1\|_{K_0 \cap \Om{2}}^2 \bigg) + \sum_{K_G \in \mathcal{T}_{G,\Gamma}} \|\nab v\|_{K_G}^2 
\end{split} 
\label{invineqgamlemfin}
\end{equation}
%
The resulting terms may be estimated by the right-hand side of \eqref{invineqgamlemres}.
\end{proof}
\end{lemma}

\subsection{Temporal estimates}

\begin{figure}[h]
\centering
\def\svgwidth{0.75\textwidth}
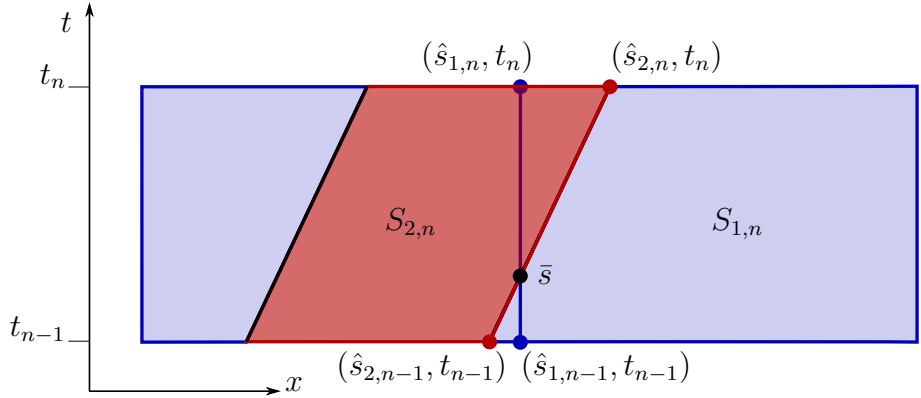
\caption{Slabwise space-time trajectories through a point $\bs \in \bGn$ for $d = 1$.}
\label{fig_lemma_disctemp_invest}
\end{figure}

Recall the domain-dependent velocity $\mu_i$, defined by \eqref{def:mui}. For a time $t^* \in I_n$, a point $x \in \Om{i}(t^*)$ and a point $s \in \Gamma(t^*)$, approached from $\Om{i}(t^*)$, we define the spatial components $\hat{x}(t)$ and $\hat{s}_i(t)$ of the slabwise space-time trajectory through $x$ and that through $s$, respectively, by
\begin{align}
\hat{x}(t) & := x + \int_{t^*}^{t}\mu_i(\tau) \ud \tau \quad \forall t \in I_n
\label{def:xhat} \\
\hat{s}_i(t) & := s + \int_{t^*}^{t}\mu_i(\tau) \ud \tau \quad \forall t \in I_n
\label{def:shat}
\end{align}
For $i = 1$, we get a straight space-time trajectory parallel to the time axis. For $i = 2$, we simply follow a point along the space-time surface $\bGn$. See Figure~\ref{fig_lemma_disctemp_invest} for an illustration. To lighten the notation, we omit the index $i$ and the time dependence when there is no risk of confusion. Thus $(\hat{s}, t) = (\hat{s}_i(t), t)$ and $\hat{s}_k = \hat{s}_{i,k} = \hat{s}_i(t_k)$, if not explicitly stated otherwise.

\begin{lemma}[Discrete temporal inverse estimates in $\| \cdot \|_{\Omega(t)}$] \label{lem:invest_temp_l2norm}
Let $k_n$ be the length of interval $I_n$ and the scaled differential operator $D_t$ be defined by \eqref{def:Dt}. For $v \in V_h^n$, let $w = w_r = D_x^r v$, where $0 \leq r \leq p$. Then, for any $v \in V_h^n$, we have that
\begin{align}
\int_{I_n} k_n^2 \| D_t w \|_{\Omt{1} \cup \Omt{O}}^2 \ud t & \lesssim \int_{I_n}  \| w \|_{\Omt{1} \cup \Omt{O}}^2 \ud t
\label{lemres:invest_temp_l2norm_1} \\
\int_{I_n} k_n^2 \| D_t w \|_{\Omt{2}}^2 \ud t & \lesssim \int_{I_n} \| w \|_{\Omt{2}}^2 \ud t
\label{lemres:invest_temp_l2norm_2}
\end{align}
\begin{proof}
The estimates follow from applying a standard one-dimensional inverse estimate for polynomials along the space-time trajectories. The presence of $D_t$ in the $I_n$-integrals gives the correct scaling for going to the space-time trajectories and back.
\end{proof}
\end{lemma}

\begin{lemma}[An inequality for $W^{1,1}((a, b))$]\label{lem:inveq_W11}
For an open interval $(a, b)$, a point $c \in (a, b)$, and for any function $w \in W^{1,1}((a, b))$ it holds that
\begin{equation}
(b - a)w(c) \leq \int_a^b w(x) \ud x + (b - a)\int_a^b |w'(x)| \ud x
\label{lemres:inveq_W11}
\end{equation}

\begin{proof}
Consider an open interval $(\alpha, \beta) \subset (a, b)$. For an arbitrary point $y \in (\alpha, \beta)$, we use integration by parts to get
\begin{equation}
\begin{split}
(\beta - y)w(\beta^-) + (y - \alpha)w(\alpha^+) \leq \int_\alpha^\beta w(x) \ud x + (\beta - \alpha) \int_\alpha^\beta | w'(x) | \ud x
\end{split}
\label{inveq_W11_0}
\end{equation}
The left-hand side of \eqref{lemres:inveq_W11} is
\begin{equation}
\begin{split}
(b - a)w(c) = & \; (b - c)w(c) + (c - a)w(c) \\
\leq & \; \int_c^b w(x) \ud x + (b - c) \int_c^b | w'(x) | \ud x \\
& + \int_a^c w(x) \ud x + (c - a) \int_a^c | w'(x) | \ud x \\
\leq & \; \int_a^b w(x) \ud x + (b - a)\int_a^b | w'(x) | \ud x
\end{split}
\end{equation}
where we have used \eqref{inveq_W11_0} with $y = \beta^- = \beta = b$ and $\alpha = c$, and \eqref{inveq_W11_0} with $\beta = c$ and $y = \alpha^+ = \alpha = a$ to obtain the first inequality. This concludes the proof.
\end{proof}
\end{lemma}

\begin{lemma}[A discrete temporal inverse estimate in $\norma{\cdot}$] \label{lem:invest_temp_anorm}
Let $\norma{\cdot}$ be defined by \eqref{def:anorm}, $k_n$ be the length of interval $I_n$, and the scaled differential operator $D_t$ be defined by \eqref{def:Dt}. Then we have that
\begin{equation}
\int_{I_n} \norma{k_n D_t v}^2 \ud t \lesssim \int_{I_n} \norma{v}^2 \ud t \quad \forall v \in V_h^n
\label{lemres:invest_temp_anorm}
\end{equation}
\begin{proof}
We expand the left-hand side of \eqref{lemres:invest_temp_anorm} by using \eqref{def:anorm}
\begin{equation}
\begin{split}
\int_{I_n} \norma{k_n D_t v}^2 \ud t = & \underbrace{\sum_{i = 1}^2 \int_{I_n} k_n^2 \| \nab D_t v \|_{\Omt{i}}^2 \ud t}_\text{= I} + \underbrace{\int_{I_n} k_n^2 \abm \|\langle \partial_{\bn^x} D_t v \rangle \|_{-1/2,h,\Gamma(t)}^2 \ud t}_\text{= II} \\
& + \; \underbrace{\int_{I_n} k_n^2 \abm \|[ D_t v ] \|_{1/2,h,\Gamma(t)}^2 \ud t}_\text{= III} + \underbrace{\int_{I_n} k_n^2 \|[\nab D_t v]\|_{\Omt{O}}^2 \ud t}_\text{= IV}
\end{split}
\label{invest_temp_anorm0}
\end{equation}
We treat the terms separately, starting with the first. Using that $\nab D_t v = D_t \nab v$ and Lemma~\ref{lem:invest_temp_l2norm}, the first term in \eqref{invest_temp_anorm0} is
\begin{equation}
\begin{split}
\text{I} & = \sum_{i = 1}^2 \int_{I_n} k_n^2 \| D_t \nab v \|_{\Omt{i}}^2 \ud t \lesssim \int_{I_n} \sum_{i=1}^2 \| \nab v \|_{\Omt{i}}^2 + \| [\nab v] \|_{\Omt{O}}^2 \ud t \\
& \leq \int_{I_n} \norma{v}^2 \ud t
\end{split}
\label{invest_temp_anorm0_I}
\end{equation}
The second term in \eqref{invest_temp_anorm0} receives the same treatment after first using Lemma~\ref{invineqgamlem}, thus
\begin{equation}
\begin{split}
\text{II} \lesssim \int_{I_n} k_n^2 \abm \bigg(\sum_{i=1}^2 \|\nab D_t v \|_{\Omt{i}}^2 + \|[\nab D_t v]\|_{\Omt{O}}^2 \bigg) \ud t \lesssim \int_{I_n} \norma{v}^2 \ud t
\end{split}
\label{invest_temp_anorm0_II}
\end{equation}
The third term in \eqref{invest_temp_anorm0} requires some more work than the others. Recall the slabwise space-time trajectories through a point $\bs \in \bGn$, whose spatial components $\hat{s} = \hat{s}_i(t)$ are defined by (\ref{def:shat}). Let $S_n^q$ denote the set of points in $I_n$ corresponding to temporal degrees of freedom for $V_h^n$. Thus $S_n^0 = \{t_n\}$ and $S_n^1 = \{t_n, t_{n-1}^+\}$. For $q > 1$, interior points of $I_n$ are also included in $S_n^q$. We consider the temporal basis functions $\lambda_k \in \mathcal{P}^q(I_n)$, where every $\lambda_k$ corresponds to a point $t_k \in S_n^q$. Writing $\hat{x}_k = \hat{x}(t_k)$, where $\hat{x}$ is defined by \eqref{def:xhat}, and using a somewhat relaxed notation, any $v \in V_h^n$ may be represented as $v(x,t) = \sum_{t_k \in S_n^q} v(\hat{x}_k, t_k)\lambda_k(t)$. With simple continuous mesh motion, $\mu$ is constant along every slabwise space-time trajectory, which means that $D_t v(\hat{x}_k, t_k) = 0$. Using this together with the somewhat relaxed representation, we have that
\begin{equation}
D_t v(x,t) = \sum_{t_k \in S_n^q} v(\hat{x}_k, t_k) D_t \lambda_k(t) = \sum_{t_k \in S_n^q} v(\hat{x}_k, t_k) \lambda_k'(t)
\label{vVhnrep_temp_Dt}
\end{equation}
With \eqref{vVhnrep_temp_Dt}, the third term in (\ref{invest_temp_anorm0}) is
\begin{equation}
\begin{split}
\text{III} & \leq \frac{k_n^2}{h_\text{min}} \int_{I_n} \abm \int_{\Gamma(t)} | (D_t v(s, t))_1 - (D_t v(s, t))_2 |^2 \ud s \ud t \\
& \leq \frac{k_n^2}{h_\text{min}} \int_{I_n} \abm \int_{\Gamma(t)} \bigg( \sum_{t_k \in S_n^q} | v_1(\hat{s}_{1, k}, t_k) - v_2(\hat{s}_{2, k}, t_k)| \underbrace{|\lambda_k'(t)|}_{\leq C(q)/k_n} \bigg)^2 \ud s \ud t \\
& \lesssim h_\text{min}^{-1} \sum_{t_k \in S_n^q} \underbrace{\int_{I_n} \abm \int_{\Gamma(t)} | v_1(\hat{s}_{1, k}, t_k) - v_2(\hat{s}_{2, k}, t_k)|^2 \ud s \ud t}_{= \text{III}.k} \\
\end{split}
\label{invest_temp_anorm0_III0}
\end{equation}
We split $\text{III}.k$ by 
\begin{equation}
\begin{split}
 \text{III}.k \lesssim & \; \underbrace{\int_{I_n} \abm \int_{\Gamma(t)} | v_1(\hat{s}_{1, k}, t_k) - v_1(\hat{s}_{2, k}, t_k)|^2 \ud s \ud t}_{=  \text{III}.k.1} \\
& + \underbrace{\int_{I_n} \abm \int_{\Gamma(t)} |v_1(\hat{s}_{2, k}, t_k) - v_2(\hat{s}_{2, k}, t_k)|^2 \ud s \ud t}_{= \text{III}.k.2} \\
\end{split}
\label{invest_temp_anorm0_IIIk}
\end{equation}
For the first term in \eqref{invest_temp_anorm0_IIIk}, we consider the spatial plane curve resulting from projecting $s(\tau) \in \Gamma(\tau)$, for all $\tau$ between $t_k$ and $t$, onto the spatial plane at time $t_k$. By applying the fundamental theorem of calculus for line integrals to this curve, we have that
\begin{equation}
\begin{split}
\text{III}.k.1 & = \int_{I_n} \abm \int_{\Gamma(t)} | v_1(s(t), t_k) - v_1(s(t_k), t_k)|^2 \ud s \ud t \\
& = \int_{I_n} \abm \int_{\Gamma(t)} \bigg| \int_{t_k}^{t} \mu(\tau) \cdot \nab v_1(s(\tau), t_k) \ud \tau \bigg|^2 \ud s \ud t \\
& \leq | \mu |_{I_n}^2 k_n \int_{I_n} \abm \int_{\Gamma(t)} \int_{I_n} | \nab v_1(s(\tau), t_k)|^2 \ud \tau \ud s \ud t \\
& = | \mu |_{I_n}^2 k_n^2 \abm \int_{\Gamma(t_k)} \int_{I_n} | \nab v_1(s(\tau), t_k)|^2 \ud \tau \ud s \\
& \lesssim | \mu |_{I_n}^2 k_n^2 \abm \int_{\Om{1}(t_k) \cup \Om{O}(t_k)} | \nab v_1(x, t_k)|^2 \ud x \\
& \lesssim | \mu |_{I_n}^2 k_n \int_{I_n} \| \nab v_1 \|_{\Om{1}(t) \cup \Om{O}(t)}^2 \ud t \\
& \lesssim | \mu |_{I_n}^2 k_n \int_{I_n} \norma{v}^2 \ud t
\end{split}
\label{invest_temp_anorm0_IIIk1}
\end{equation}
where, in the fifth step, we have taken possible multiples of the same line integrals into account and expanded the domain of integration. In the sixth step, we have used a standard inverse inequality for polynomials. For the second term in \eqref{invest_temp_anorm0_IIIk}, we use Lemma~\ref{lem:inveq_W11}, thus
\begin{equation}
\begin{split}
\text{III}.k.2 & = k_n \abm \int_{\Gamma(t_k)} | v_1(s, t_k) - v_2(s, t_k)|^2 \ud s = \int_{\Gamma(t_k)} \bigg(k_n \abm [v]^2(s, t_k)\bigg) \ud s \\
& \lesssim \int_{\Gamma(t_k)} \bigg(\int_{I_n} \abm [v]^2(\hat{s}_2, t) \ud t + k_n \abm \int_{I_n} \bigg| D_{t,2}[v]^2(\hat{s}_2, t) \bigg| \ud t \bigg) \ud s \\
& \lesssim \int_{\bGn} |[v]|^2 \ud \bs + k_n \int_{\bGn} | [v] | |D_{t,2} [v] | \ud \bs \\
& \leq \| [v] \|_{\bGn}^2 + k_n \| [v] \|_{\bGn} \| D_{t,2} [v] \|_{\bGn} \\
& \leq \underbrace{\bigg(1 + \frac{1}{\varepsilon} \bigg) \| [v] \|_{\bGn}^2}_{= \text{III}.k.2.1} + \underbrace{\varepsilon k_n^2 \| D_{t,2} [v] \|_{\bGn}^2}_{= \text{III}.k.2.2}
\end{split}
\label{invest_temp_anorm0_IIIk2}
\end{equation}
Using \eqref{HHnormineqs}, the first term is
\begin{equation}
\begin{split}
\text{III}.k.2.1 & = \bigg(1 + \frac{1}{\varepsilon} \bigg) \int_{I_n} \abm \| [v] \|_{\Gamma(t)}^2 \ud t \lesssim h \int_{I_n} \abm \| [v] \|_{1/2, h, \Gamma(t)}^2 \ud t \\
& \leq h \int_{I_n} \norma{v}^2 \ud t
\end{split}
\label{invest_temp_anorm0_IIIk2}
\end{equation}
Using again \eqref{HHnormineqs}, the second term is
\begin{equation}
\begin{split}
 \text{III}.k.2.2 & = \varepsilon k_n^2 \| D_{t,2} v_1 - D_{t,2} v_2 \|_{\bGn}^2 = \varepsilon k_n^2 \| D_t v_1 + \mu \cdot \nab v_1 - D_t v_2 \|_{\bGn}^2 \\
& \lesssim \varepsilon k_n^2 \| [D_t v] \|_{\bGn}^2 + \varepsilon k_n^2 \| \mu \cdot \nab v_1 \|_{\bGn}^2 \\
& \leq \varepsilon k_n^2 \int_{I_n} \abm \| [D_t v] \|_{\Gamma(t)}^2 \ud t + |\mu|_{I_n}^2 k_n^2 \| \nab v_1 \|_{\bGn}^2 \\
& \leq \varepsilon h \underbrace{\int_{I_n} k_n^2 \abm \| [D_t v] \|_{1/2, h, \Gamma(t)}^2 \ud t}_{= \text{III}} + \underbrace{|\mu|_{I_n}^2 k_n^2 \| \nab v_1 \|_{\bGn}^2}_{= \text{III}.k.2.2.2}
\end{split}
\label{invest_temp_anorm0_IIIk22}
\end{equation}
where the first term is done. For the second term, we use Corollary~\ref{cor:scatraineqGamK}, thus
\begin{equation}
\begin{split}
\text{III}.k.2.2.2 & = |\mu|_{I_n}^2 k_n^2 \int_{I_n} \abm \sum_{K \in \mathcal{T}_{0, \Gamma(t)}} \| \nab v_1 \|_{\Gamma_K}^2 \ud t \\
& \lesssim \frac{|\mu|_{I_n}^2 k_n^2}{h_\text{min}} \int_{I_n} \sum_{K \in \mathcal{T}_{0, \Gamma(t)}} \| \nab v_1 \|_{K}^2 \ud t \\
& \leq \frac{|\mu|_{I_n}^2 k_n^2}{h_\text{min}} \int_{I_n} \| \nab v_1 \|_{\Omt{1} \cup \Omt{O}}^2 \ud t \lesssim \frac{|\mu|_{I_n}^2 k_n^2}{h_\text{min}}  \int_{I_n} \norma{v}^2 \ud t
\end{split}
\label{invest_temp_anorm0_IIIk222}
\end{equation}
This concludes the separate treatment of all the terms unfolding in the estimation of the third term in \eqref{invest_temp_anorm0}. Collecting all the estimates and using \eqref{quasiuniformity_st} gives us
\begin{equation}
\begin{split}
\text{III} & \lesssim h_\text{min}^{-1} \sum_{t_k \in S_n^q} \bigg( | \mu |_{I_n}^2 k_n \int_{I_n} \norma{v}^2 \ud t + h \int_{I_n} \norma{v}^2 \ud t + \text{III}.k.2.2  \bigg) \\
& \lesssim h_\text{min}^{-1} \bigg( \bigg(| \mu |_{I_n}^2 k_n + h + \frac{|\mu|_{I_n}^2 k_n^2}{h_\text{min}} \bigg)\int_{I_n} \norma{v}^2 \ud t + \varepsilon h (\text{III})\bigg) \\
& \lesssim \int_{I_n} \norma{v}^2 \ud t + \varepsilon (\text{III})
\end{split}
\label{invest_temp_anorm0_III_fin0}
\end{equation}
By kicking back the $\varepsilon$-term and taking $\varepsilon$ sufficiently small, we may estimate the third term in \eqref{invest_temp_anorm0} by the first term on the right-hand side of \ref{invest_temp_anorm0_III_fin0}. The fourth term in \eqref{invest_temp_anorm0} receives the same treatment as the first, thus
\begin{equation}
\begin{split}
\text{IV} \lesssim \int_{I_n} k_n^2 \bigg( \| D_t \nab v \|_{\Omt{1} \cup \Omt{O}}^2 + \| D_t \nab v \|_{\Omt{2}}^2 \bigg) \ud t \lesssim \int_{I_n} \norma{v}^2 \ud t
\end{split}
\label{invest_temp_anorm0_IV}
\end{equation}
The treatment of all the terms in \eqref{invest_temp_anorm0} is done. This shows \eqref{lemres:invest_temp_anorm}.
\end{proof}
\end{lemma}

\begin{lemma}[An inverse inequality for $\mathcal{P}((a, c), (c, b))$] \label{lem:inveq_poly_disc}
For an open interval $(a, b)$, a point $c \in (a, b)$, and for $w \in \mathcal{P}((a, c), (c, b))$, i.e., $w$ is a polynomial on $(a, c)$, possibly another polynomial on $(c, b)$, and possibly discontinuous at $c$, there exists a positive constant depending on the polynomial degree such that
\begin{equation}
\begin{split}
(b - a)|w(a^+)|^2 \lesssim & \; \int_a^c |w(x)|^2 \ud x + \int_c^b |w(x)|^2 \ud x + (b - c)|[w](c)|^2 \\
& + (b - c)(c - a)\int_a^c |w'(x)|^2 \ud x
\label{lemres:inveq_poly_disc}
\end{split}
\end{equation}

\begin{proof}
Using a standard inverse inequality for polynomials with a positive constant that depends on the polynomial degree, the left-hand side of \eqref{lemres:inveq_poly_disc} is
\begin{equation}
\begin{split}
(b - a)|w(a^+)|^2 & \lesssim \int_a^c |w(x)|^2 \ud x + (b - c)|w(a^+)|^2
\label{inveq_poly_disc0}
\end{split}
\end{equation}
Adding and subtracting $w(c^-)$ and $w(c^+)$ within the absolute value, followed by using standard estimates, the second term is
\begin{equation}
\begin{split}
(b - c)|w(a^+)|^2 = & \; (b - c)\bigg| -\int_a^c w'(x) \ud x - [w](c) + w(c^+)\bigg|^2 \\
\lesssim & \; (b - c)(c - a)\int_a^c |w'(x)|^2 \ud x \\
& + (b - c)|[w](c)|^2 + \int_c^b |w(x)|^2 \ud x
\end{split}
\label{inveq_poly_disc0_II}
\end{equation}
\end{proof}
\end{lemma}

\begin{lemma}[Discrete temporal inverse inequalities for $V_h^n$] \label{lem:inv_ineq_temp_om}
Let $k_n$ be the length of interval $I_n$, the scaled differential operator $D_t$ be defined by \eqref{def:Dt}, $|\mu|_{I_n} = \max_{t \in I_n} \{|\mu(t)|\}$, and $\norma{\cdot}$ be defined by \eqref{def:anorm}. Then, for $q = 0, 1$, we have that
\begin{align}
\begin{split}
\sum_{i=1}^2 k_{n}^2 \| (D_t v)_{n-1}^+\|_{\Om{{i, n-1}}}^2 \lesssim & \; \sum_{i=1}^2 \int_{I_n} k_n \| D_t v\|_{\Omt{i}}^2 \ud t \\
& + |\mu|_{I_n} \int_{I_n} \norma{v}^2 \ud t \quad \forall v \in V_h^n
\end{split}
\label{lemres:inv_ineq_temp_omp} \\
\begin{split}
\sum_{i=1}^2 k_{n}^2 \| (D_t v)_{n}^-\|_{\Om{{i,n}}}^2
\lesssim & \; \sum_{i=1}^2 \int_{I_n} k_n \| D_t v\|_{\Omt{i}}^2 \ud t \\
& + |\mu|_{I_n} \int_{I_n} \norma{v}^2 \ud t \quad \forall v \in V_h^n
\end{split}
\label{lemres:inv_ineq_temp_omm} 
\end{align}
\begin{proof}
We only prove (\ref{lemres:inv_ineq_temp_omp}), since the proof of (\ref{lemres:inv_ineq_temp_omm}) is analogous. Recall $\hat{x}(t)$ defined by \eqref{def:xhat}. We denote by $\bar{\hat{x}}_i^n$ the slabwise space-time trajectory  through a point $x_{i, n-1} \in \Om{{i, n-1}}$. We define the set of points in $\Om{{1, n-1}}$ with cut and uncut space-time trajectories by
\begin{align}
\Omega_{1, n-1}^{\Gamma} & := \{ x \in \Om{{1, n-1}} : \bar{\hat{x}}_1^n \cap \bGn \neq \emptyset \}
\label{def:setOm1nm1G} \\
\Omega_{1, n-1}^{n} & := \{ x \in \Om{{1, n-1}} : \bar{\hat{x}}_1^n \cap \bGn = \emptyset \}
\label{def:setOm1nm1n}
\end{align}
The idea to prove \eqref{lemres:inv_ineq_temp_omp} is that if a point's space-time trajectory is \emph{uncut}, we use a standard inverse inequality, and if it is \emph{cut}, we use Lemma~\ref{lem:inveq_poly_disc}. Using that $\Omega_{1, n-1}^{\Gamma}$ and $\Omega_{1, n-1}^{n}$ form a partition of $\Om{{1, n-1}}$, the left-hand side of (\ref{lemres:inv_ineq_temp_omp}) is
\begin{equation}
\begin{split}
\sum_{i=1}^2 k_{n}^2 \| (D_t v)_{n-1}^+\|_{\Om{{i,n-1}}}^2 = & \; \underbrace{k_{n}^2 \| (D_t v)_{n-1}^+\|_{\Omega_{1, n-1}^{n}}^2}_{= \text{I}} + \underbrace{k_{n}^2 \| (D_t v)_{n-1}^+\|_{\Om{{2,n-1}}}^2}_{= \text{II}} \\
& + \underbrace{k_{n}^2 \| (D_t v)_{n-1}^+\|_{\Omega_{1, n-1}^{\Gamma}}^2}_{= \text{III}}
\end{split}
\label{inv_ineq_temp_om0} 
\end{equation}
We treat the terms separately. See Figure~\ref{fig_lemma_disctemp_invineq} for an illustration of the proof idea.
\begin{figure}[h]
\centering
\def\svgwidth{0.85\textwidth}
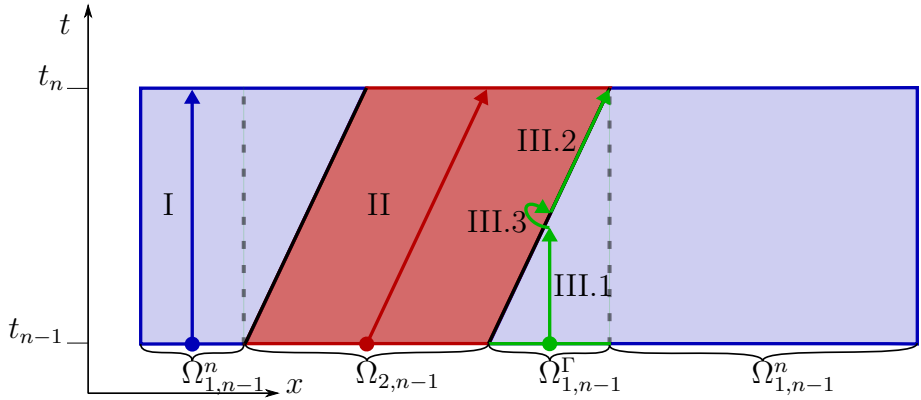
\caption{The starting domains $\Omega_{1, n-1}^{n}$, $\Omega_{2, n-1}$, and $\Omega_{1, n-1}^{\Gamma}$. The arrows represent the treatment of the corresponding right-hand side terms.}
\label{fig_lemma_disctemp_invineq}
\end{figure}
\noindent Using a standard inverse inequality, the first and second term in \eqref{inv_ineq_temp_om0} are
\begin{equation}
\text{I} \lesssim \int_{I_n} k_n \| D_t v \|_{\Omt{1}}^2 \ud t \qquad \text{II} \lesssim \int_{I_n} k_n \| D_t v \|_{\Omt{2}}^2 \ud t
\label{inv_ineq_temp_om0_IandII} 
\end{equation}
For the third term in \eqref{inv_ineq_temp_om0}, we recall $\hat{s}$ defined by \eqref{def:shat}. We consider a space-time curve that starts at $x \in \Omega_{1, n-1}^{\Gamma}$, goes straight up in time until it hits $\bGn$, which occurs at time $t_\Gamma$, then travels on $\bGn$ along $(\hat{s}(t), t)$ up to $t_n$. We will apply Lemma~\ref{lem:inveq_poly_disc} to the function that is $(D_t v)_1$ up until $t_\Gamma$ along this space-time curve, and $(D_t v)_2$ afterwards. Here the corresponding derivative term on the right-hand side of \eqref{lemres:inveq_poly_disc} vanishes since $D_t^2 v(x, t) = 0$ for $q \leq 1$. Thus
\begin{equation}
\begin{split}
\text{III} = & \; k_{n} \int_{\Omega_{1, n-1}^{\Gamma}} \bigg( k_n |D_t v(x, t_{n-1}^+)|^2 \ud x \bigg) \\
\lesssim & \; k_{n} \int_{\Omega_{1, n-1}^{\Gamma}} \bigg( \int_{t_{n-1}}^{t_{\Gamma}} |D_t v(x, t)|^2 \ud t + \int_{t_{\Gamma}}^{t_n} \abm |(D_t v)_2(\hat{s}(t), t)|^2 \ud t \\
& + (t_n - t_{\Gamma})|[D_t v](x, t_{\Gamma})|^2 + (t_n - t_{\Gamma})(t_{\Gamma} - t_{n-1})\int_{t_{n-1}}^{t_{\Gamma}} | D_t^2 v(x, t) |^2 \ud t \bigg) \ud x \\
\lesssim & \; \underbrace{k_{n} \int_{\Omega_{1, n-1}^{\Gamma}} \int_{t_{n-1}}^{t_{\Gamma}} |D_t v(x, t)|^2 \ud t \ud x}_{= \text{III.1}} + \underbrace{k_{n}\int_{\Omega_{1, n-1}^{\Gamma}} \int_{t_{\Gamma}}^{t_n} \abm |(D_t v)_2(\hat{s}(t), t)|^2 \ud t \ud x}_{= \text{III.2}} \\
& + \underbrace{k_{n} \int_{\Omega_{1, n-1}^{\Gamma}} (t_n - t_{\Gamma})|[D_t v](x, t_{\Gamma})|^2 \ud x}_{= \text{III.3}}.
\end{split}
\label{inv_ineq_temp_om0_III} 
\end{equation}
Simply expanding the domain of integration, the first term is
\begin{equation}
\text{III}.1 \lesssim \int_{I_n} k_n \| D_t v \|_{\Omt{1}}^2 \ud t
\label{inv_ineq_temp_om0_III1} 
\end{equation}
For the second and third term in \eqref{inv_ineq_temp_om0_III}, we want to change the domain of integration from $\Omega_{1, n-1}^{\Gamma}$ to its temporal projection onto $\bGn$. To do this, we note that $\ud x \lesssim |\mu|_{I_n} \ud \bs$, where $\ud x$ and $\ud \bs$ are the integration differentials for $\Omega_{1, n-1}^{\Gamma}$ and $\bGn$, respectively. Using this, a standard trace inequality for $H^1(\Om{2})$, that $\nab D_t v = D_t \nab v$, and Lemma~\ref{lem:invest_temp_l2norm}, the second term is
\begin{equation}
\begin{split}
\text{III}.2 & \lesssim k_{n} \int_{\Omega_{1, n-1}^{\Gamma}} \int_{I_n} \abm |(D_t v)_2(\hat{s}(t), t)|^2 \ud t \ud x \\
& \lesssim k_{n} \int_{\bGn} \int_{I_n} \abm |(D_t v)_2(\bs)|^2 \ud t |\mu|_{I_n} \ud \bs \lesssim |\mu|_{I_n} k_{n}^2 \int_{I_n} \abm \| (D_t v)_2 \|_{\Gamma(t)}^2 \ud t \\
& \lesssim |\mu|_{I_n} k_{n} \int_{I_n} k_n \| D_t v \|_{\Omt{2}}^2 \ud t + |\mu|_{I_n} \int_{I_n} k_{n}^2 \| D_t \nab v \|_{\Omt{2}}^2 \ud t \\
& \lesssim \int_{I_n} k_n \| D_t v \|_{\Omt{2}}^2 \ud t + |\mu|_{I_n} \int_{I_n} \norma{v}^2 \ud t
\end{split}
\label{inv_ineq_temp_om0_III2} 
\end{equation}
Using the relation between the integration differentials, the estimate \eqref{HHnormineqs}, and Lemma~\ref{lem:invest_temp_anorm}, the third term in (\ref{inv_ineq_temp_om0_III}) is
\begin{equation}
\begin{split}
\text{III}.3 & \lesssim k_{n}^2 \int_{\Omega_{1, n-1}^{\Gamma}} |[D_t v](x, t_{\Gamma})|^2 \ud x \lesssim k_{n}^2 \int_{\bGn} |[D_t v](\bs)|^2 |\mu|_{I_n} \ud \bs \\
& \lesssim |\mu|_{I_n} k_{n}^2 \int_{I_n} \abm \| [D_t v] \|_{\Gamma(t)}^2 \ud t \lesssim |\mu|_{I_n} h \int_{I_n} \norma{k_{n} D_t v}^2 \ud t \\
& \lesssim |\mu|_{I_n} \int_{I_n} \norma{v}^2 \ud t
\end{split}
\label{inv_ineq_temp_om0_III3} 
\end{equation}
The treatment of all the terms in \eqref{inv_ineq_temp_om0} is done. This shows \eqref{lemres:inv_ineq_temp_omp}.
\end{proof}
\end{lemma}

\begin{lemma}[A discrete temporal inverse estimate in $\normx{\cdot}$] \label{lem:invest_temp_xnorm}
Let the norm $\normx{\cdot}$ be defined by \eqref{def:xnorm}, $k_n$ be the length of time interval $I_n$, and the scaled differential operator $D_t$ be defined by \eqref{def:Dt}. Then, for $q = 0, 1$, we have that
\begin{equation}
\normx{k_n D_t v} \lesssim \normx{v} \quad \forall v \in V_h
\label{lemres:invest_temp_xnorm}
\end{equation}

\begin{proof}
The square of the left-hand side of \eqref{lemres:invest_temp_xnorm} is
\begin{equation}
\begin{split}
\normx{k_n D_t v}^2 = & \; \sum_{i=1}^2 \sum_{n=1}^N \int_{I_n} k_n\| D_t (k_n D_t v) \|_{\Omt{i}}^2 \ud t \\
& + \sum_{n=1}^N \bigg(\int_{I_n} \norma{k_n D_t v}^2 \ud t + \| |\bn^t|^{1/2} [k_n D_t v]\|_{\bGn}^2 \bigg) \\
& + \sum_{i=1}^2 \sum_{n=1}^{N-1} \|[k_n D_t v]_n \|_{\Om{{i,n}}}^2 \\
& + \sum_{i=1}^2 \bigg( \| (k_n D_t v)_N^- \|_{\Om{{i,N}}}^2 + \| (k_n D_t v)_0^+ \|_{\Om{{i,0}}}^2 \bigg)
\end{split}  
\label{invest_temp_xnorm0} 
\end{equation}
The first term vanishes since $D_t^2 v(x, t) = 0$ for $q \leq 1$. The $\bGn$-norm term is estimated by the $A_{h,t}$-norm term by using \eqref{HHnormineqs}. Applying Lemma~\ref{lem:invest_temp_anorm} to the $A_{h,t}$-norm term and Lemma~\ref{lem:inv_ineq_temp_om} to all the terms in the last two rows, we get terms which may be estimated by $\normx{v}^2$.
\end{proof}
\end{lemma}

\section{Interpolation} \label{sec:interpolation}

Let $\cdot^\circ$ denote the interior of a set, e.g., $I_n^\circ = (t_{n-1}, t_n)$. Also let $C_b(\cup_n I_n^\circ)$ denote the space of functions that are continuous and bounded on every $I_n^\circ$. In this section, the space-time interpolation operator $\interph : C_b(\cup_n I_n^\circ; L^1(\Om{0})) \to V_h$ is successively constructed by first defining spatial interpolation operators, then temporal ones, and finally combining them. Interpolation error estimates are also presented.

\subsection{Slabwise operators and local estimates}

\begin{definition}(Spatial interpolation operators)
\label{def:spatial_interpop}
We define the spatial interpolation operators $\pi_{h,0} : L^1(\Om{0}) \to V_{h,0}$ and $\pi_{h,G} : L^1(G) \to V_{h,G}$ to be the Scott-Zhang interpolation operators for the spaces $V_{h,0}$ and $V_{h,G}$, respectively, where the defining integrals are taken over entire simplices.
\end{definition}
\noindent Note that $\pi_{h,G}$ is time-dependent but to lighten the notation we omit this. The temporal interpolation operators will interpolate along the space-time trajectories of the domains $\Om{0}$ and $G$. For $n = 1, \dots, N$, we define the slabwise space-time trajectory for a point $x \in \Om{0}$ and that of a point $x_n \in G(t_n)$ by
\begin{align}
\bar{\hat{x}}_0^n & := \{(\hat{x}(t), t) : \hat{x}(t) = x, t \in I_n\}
\label{def:st_traject_0} \\
\bar{\hat{x}}_G^n & := \{(\hat{x}(t), t) : \hat{x}(t) = x_n - \int_{t}^{t_n} \mu(\tau) \ud \tau, t \in I_n\}
\label{def:st_traject_G} 
\end{align}
Note that \eqref{def:xhat} can be used to obtain all trajectories defined by \eqref{def:st_traject_G} but \emph{not} all defined by \eqref{def:st_traject_0} because some $\bar{\hat{x}}_0^n$ may lie completely in $S_{2,n}$. Let $S_n^q$ denote the set of temporal interpolation points for interpolation to $\mathcal{P}^q(I_n)$. We take $S_n^0 = \{t_n^-\}$ and $S_n^1 = \{t_n^-, t_{n-1}^+\}$. For $q > 1$, we include interior points of $I_n$ in some suitable fashion.
\begin{definition}(Temporal interpolation operators)
\label{def:temporal_interpop}
For each time subinterval $I_n$, where $n = 1, \dots, N$, we define the temporal interpolation operators $\pi_{0}^n : C_b({\bar{\hat{x}}_0^n}^\circ) \to \mathcal{P}^q(\bar{\hat{x}}_0^n)$ and $\pi_{G}^n : C_b({\bar{\hat{x}}_G^n}^\circ) \to \mathcal{P}^q(\bar{\hat{x}}_G^n)$ to be the nodal interpolation operators that use the points in $S_n^q$ as nodal interpolation points.
\end{definition}
\noindent 
Note that $\pi_{0}^n$ and $\pi_{G}^n$ are spatially dependent but to lighten the notation we omit this. We combine the spatial and temporal interpolation operators to define space-time ones.
\begin{definition}(Slabwise space-time interpolation operators)
For $n = 1, \dots, N$, we define the slabwise space-time interpolation operators $\interp{{h,0}}{n} : C_b(I_n^\circ; L^1(\Om{0})) \to V_{h,0}^n$ and $ \interp{{h,G}}{n} : C_b(I_n^\circ; L^1(G)) \to V_{h,G}^n$ by
\begin{equation}
\interp{{h,0}}{n} := \pi_{0}^n \pi_{h,0}  \qquad \interp{{h,G}}{n} := \pi_{G}^n \pi_{h,G}
\label{def_interphn}
\end{equation} 
\end{definition}
\noindent Recall the interdependent indices $i \in \{1, 2\}$ and $j \in \{0, G\}$ where $j = 0$ for $ i=1$ and $j = G$ for $i=2$. Let $\bar{K}_n := \{(x, t) : x \in K=K(t), t \in I_n\}$ denote an arbitrary space-time prism, where $K = K_j \in \mathcal{T}_j$. Let $\mathcal{N}(K)$ denote the neighborhood of a simplex $K$, i.e., the set of all adjacent simplices to and including $K$. We also use the notation $\| w \|_{K, I_n} = \max_{t \in I_n} \{ \| w(\cdot, t) \|_{K(t)} \}$.

\begin{lemma}[Local space-time interpolation error estimates for $\bar{K}_n$]
\label{lem:interpest_loc_st_bKn}
Let $\interp{{h,j}}{n}$ be defined by \eqref{def_interphn}, where $j \in \{0, G\} $, and let $D_t$ be defined by \eqref{def:Dt}. Then, for a function $v$ with sufficient spatial and temporal regularity, we have for $0 \leq s \leq q + 1$ and $0 \leq r \leq p + 1$ that
\begin{align}
\| D_t^s (v - \interp{{h,j}}{n} v) \|_{\bar{K}_n} \lesssim k_n^{q+1-s}\| D_t^{q+1} v \|_{\bar{K}_n} + h^{p+1} k_n^{1/2} \| D_x^{p+1} D_t^s v \|_{\mathcal{N}(K), I_n}
\label{interpest_loc_st_t_res} \\
\| D_x^r (v - \interp{{h,j}}{n} v) \|_{\bar{K}_n} \lesssim k_n^{q+1}\| D_t^{q+1} D_x^r v \|_{\bar{K}_n} + h^{p+1-r} k_n^{1/2} \| D_x^{p+1} v \|_{\mathcal{N}(K), I_n}
\label{interpest_loc_st_x_res}
\end{align}
\begin{proof}
We show the two estimates separately, starting with the first. Using that $\interp{{h,j}}{n} = \pi_{j}^n \pi_{h, j} = \pi_{h, j} \pi_{j}^n$, that $D_t^s\pi_{h,j} = \pi_{h,j}D_t^s$, stability of $\pi_{h, j}$, and a trivial estimate, the left-hand side of \eqref{interpest_loc_st_t_res} is
\begin{equation}
\begin{split}
\| D_t^s(v - \interp{{h,j}}{n} v) \|_{\bar{K}_n} & \leq \| D_t^s \pi_{h,j} ( \mathds{1} - \pi_{j}^n) v \|_{\bar{K}_n} + \| D_t^s(v - \pi_{h,j} v ) \|_{\bar{K}_n} \\
& \lesssim \| D_t^s( \mathds{1} - \pi_{j}^n) v \|_{\bar{K}_n} + k_n^{1/2} \| D_t^s v - \pi_{h,j} D_t^s v \|_{K, I_n}
\end{split}
\label{interpest_loc_st_t_0} 
\end{equation}
Applying standard estimates for $\pi_j^n$ and $\pi_{h, j}$ shows the first estimate and we move on to the second. We are going to use the expansion of interpolants of $\pi_j^n$ into a sum over the temporal interpolation points $t_k \in S_n^q$ with $\lambda_k \in \mathcal{P}^{q}(I_n)$ denoting the corresponding shape function. For a function $w$ of sufficient regularity, we have that
\begin{equation}
\begin{split}
\| \pi_j^n w \|_{\bar{K}_n}^2 & = \int_{I_n} \int_{K(t)} \bigg| \sum_{t_k \in S_n^q} w(\hat{x}(t_k), t_k)\lambda_k(t) \bigg|^2 \ud x \ud t \\
& \leq (q+1) \sum_{t_k \in S_n^q} \int_{I_n} \int_{K(t)} | w(\hat{x}(t_k), t_k) |^2 \ud x \ud t \lesssim k_n \| w \|_{K, I_n}^2
\end{split}
\label{interp_tempexp_L2ns}
\end{equation}
Using this after using that $D_x^r \pi_j^n = \pi_j^n D_x^r$, the left-hand side of \eqref{interpest_loc_st_x_res} is
\begin{equation}
\begin{split}
\| D_x^r(v - \interp{{h,j}}{n} v) \|_{\bar{K}_n} & \leq \| D_x^r(v - \pi_{j}^n v ) \|_{\bar{K}_n} + \| D_x^r \pi_j^n ( \mathds{1} - \pi_{h, j}) v \|_{\bar{K}_n} \\
& \lesssim \| D_x^rv - \pi_{j}^n D_x^r v \|_{\bar{K}_n} + k_n^{1/2} \| D_x^r( \mathds{1} - \pi_{h, j}) v \|_{K, I_n} 
\end{split}
\label{interpest_loc_st_x_0} 
\end{equation}
Applying standard estimates for $\pi_j^n$ and $\pi_{h, j}$ shows the second estimate.
\end{proof}
\end{lemma}

\noindent Recall $\hat{s}$ defined by \eqref{def:shat}. Let $\mathcal{T}_{j, D} = \{ K \in \mathcal{T}_j : K \cap D \neq 0\}$, where $D$ is a possibly time-dependent subset of $\mathbb{R}^{d+1}$.

\begin{lemma}[Slabwise space-time interpolation error estimates for $\bGn$]
\label{lem:interpest_loc_st_bGn}
Let $\interp{{h,j}}{n}$ be defined by \eqref{def_interphn}, where $j \in \{0, G\}$, and let $D_t$ be defined by \eqref{def:Dt}. Then, for any function $v$ with sufficient spatial and temporal regularity, we have that
\begin{equation}
\| (v - \interp{{h,j}}{n} v)_i \|_{\bGn}^2 \lesssim k_n^{2q+2} \| D_t^{q+1} v \|_{L^2(\bGn, L^\infty(I_n))}^2 + h^{2p+1} \sum_{K \in \mathcal{T}_{j, \bGn}} k_n \| D_x^{p+1}v \|_{\mathcal{N}(K), I_n}^2
\label{lemres:interpest_loc_st_bGn}
\end{equation} 
\begin{proof}
The general proof idea is the same for all $q \geq 0$. What varies is how a temporal difference is treated. We show how to treat if for $q=1$ from which it should be relatively straightforward how to handle the other cases. Using the shape functions $\lambda_k \in \mathcal{P}^q(I_n)$, corresponding to interpolation points $t_k \in S_n^q$, the argument of the norm on the left-hand side of \eqref{lemres:interpest_loc_st_bGn} is
\begin{equation}
\begin{split}
& \; (v - \interp{{h,j}}{n} v)_i |_{\bGn} = v(s, t) - \sum_{t_k \in S_n^q}\pi_{h,j}v(\hat{s}_k, t_k)\lambda_k(t) \\
= & \; \underbrace{\sum_{t_k \in S_n^q}\bigg(v(s, t) - v(\hat{s}_k, t_k)\bigg)\lambda_k(t)}_{= A} + \underbrace{\sum_{t_k \in S_n^q}\bigg(v(\hat{s}_k, t_k) - \pi_{h,j}v(\hat{s}_k, t_k)\bigg)\lambda_k(t)}_{= B}
\end{split}
\label{interpest_loc_st_gamma_q1_split_0} 
\end{equation}
The left-hand side of \eqref{lemres:interpest_loc_st_bGn} may thus be split by $\| (v - \interp{{h,j}}{n} v)_i \|_{\bGn}^2 \lesssim \| A \|_{\bGn}^2 + \| B \|_{\bGn}^2$ where we consider the terms separately, starting with the first. We proceed with some further treatment of $A$ for which we restrict ourselves to the case $q = 1$. From this case it should however be relatively straightforward how to treat $A$ for $q \neq 1$. Using the mean value theorem along the space-time trajectories, the explicit expressions for the shape functions $\lambda_{n-1}$ and $\lambda_n$ for $q = 1$, and the fundamental theorem of calculus, we have
\begin{equation}
A = \sum_{t_k \in S_n^q}D_tv(\hat{s}, c_k)(t - t_k)\lambda_k(t) = \frac{(t - t_{n})(t - t_{n-1})}{k_n}\int_{c_{n-1}}^{c_n} D_{t}^2v(\hat{s}, \tau) \ud \tau
\label{interpest_loc_st_gamma_q1_split_A} 
\end{equation}
Using \eqref{interpest_loc_st_gamma_q1_split_A}, we have that
\begin{equation}
\begin{split}
\| A \|_{\bGn}^2 & \leq \int_{I_n} \abm \int_{\Gamma(t)} k_n^2 (c_n - c_{n-1}) \int_{c_{n-1}}^{c_n} | D_t^2 v(\hat{s}, \tau) |^2 \ud \tau \ud s \ud t \\
& \leq k_n^4 \| D_t^2 v \|_{L^2(\bGn, L^\infty(I_n))}^2
\end{split}
\label{interpest_loc_st_gamma_q1_I} 
\end{equation}
Writing $B_k = v(\hat{s}_k, t_k) - \pi_{h,j}v(\hat{s}_k, t_k)$, using \eqref{scatraineqGamK_warmup}, and standard estimates for $\pi_{h, j}$, we have that
\begin{equation}
\begin{split}
\| B \|_{\bGn}^2 & \lesssim \int_{I_n} \abm \int_{\Gamma(t)} \sum_{t_k \in S_n^q}|v(\hat{s}_k, t_k) - \pi_{h,j}v(\hat{s}_k, t_k)|^2 \ud s \ud t \\
& \lesssim \sum_{t_k \in S_n^q} \sum_{K \in \mathcal{T}_{j, \bGn}} k_n \|B_k\|_{\Gamma_K}^2 \lesssim \sum_{t_k \in S_n^q} \sum_{K \in \mathcal{T}_{j, \bGn}} k_n \bigg(h_K^{-1} \| B_k \|_K^2 + h_K \| D_x B_k \|_K^2 \bigg) \\
& \lesssim h^{2p+1} \sum_{K \in \mathcal{T}_{j, \bGn}} k_n \| D_x^{p+1}v \|_{\mathcal{N}(K), I_n}^2
\end{split}
\label{interpest_loc_st_gamma_q1_II} 
\end{equation}
\end{proof}
\end{lemma}

\begin{lemma}[Local spatial interpolation error estimates for temporal endpoints]
\label{lem:interpest_loc_s}
Let $\interp{{h,j}}{n}$ be defined by \eqref{def_interphn}, where $j \in \{0, G\}$, let $t_k^\tau \in \{t_{n-1}^+, t_n^-\}$, and let $D_t$ be defined by \eqref{def:Dt}. Then, for any function $v$ with sufficient spatial and temporal regularity, we have for $q > 0$ that
\begin{equation}
\|(v - \interp{{h,j}}{n}v)_k^\tau \|_{K} \lesssim h^{p+1} \| D_x^{p+1}v (\cdot, t_k^\tau) \|_{\mathcal{N}(K)}
\label{interpest_loc_s_gen_res}
\end{equation}
and for $q = 0$ that
\begin{align}
\|(v - \interp{{h,j}}{n}v)_n^- \|_{K} &\lesssim h^{p+1} \| D_x^{p+1}v(\cdot, t_n^-) \|_{\mathcal{N}(K)}
\label{interpest_loc_s_q0minus_res} \\
\|(v - \interp{{h,j}}{n}v)_{n-1}^+ \|_{K} &\lesssim k_n^{1/2}\| D_t v \|_{\bar{K}_n} + h^{p+1} \| D_x^{p+1}v (\cdot, t_n^-) \|_{\mathcal{N}(K)}
\label{interpest_loc_s_q0plus_res}
\end{align}
\begin{proof}
Estimates \eqref{interpest_loc_s_gen_res} and \eqref{interpest_loc_s_q0minus_res} follow from simply using that $t_k^\tau$ is an interpolation point of $\pi_j^n$ and then a standard estimate for $\pi_{h,j}$. This does not work for \eqref{interpest_loc_s_q0plus_res}, since $t_{n-1}^+$ is not an interpolation point for $q = 0$. Instead, we integrate along the slabwise space-time trajectory of an element $x \in K$ to obtain
\begin{equation}
\begin{split}
\|v(\cdot, t_{n-1}^+) - v(\cdot, t_{n}^-)\|_K^2 & = \int_K \bigg|v(x(t_{n-1}^+), t_{n-1}^+) - v(x(t_n^-), t_{n}^-)\bigg|^2 \ud x \\
& = \int_K \bigg| \int_{I_n} D_t v(x(t), t) \ud t \bigg|^2 \ud x \leq k_n \| D_t v \|_{\bar{K}_n}^2
\end{split}
\label{interpest_loc_s_q0plus_aux}   
\end{equation} 
Using the definition of $\pi_j^n$ for $q=0$, the left-hand side of \eqref{interpest_loc_s_q0plus_res} is
\begin{equation}
\begin{split}
\|(v - \interp{{h,j}}{n}v)_{n-1}^+ \|_{K} & \leq \|(v - \pi_j^n v)(\cdot, t_{n-1}^+)\|_K + \|(\pi_j^n v - \pi_j^n\pi_{h,j}v)(\cdot, t_{n-1}^+) \|_{K} \\
& = \|v(\cdot, t_{n-1}^+) - v(\cdot, t_{n}^-)\|_K + \|(v - \pi_{h,j}v)(\cdot, t_{n}^-) \|_{K}
\end{split}
\label{interpest_loc_s_q0plus_res_0}   
\end{equation} 
Applying \eqref{interpest_loc_s_q0plus_aux} and a standard estimate for $\pi_{h,j}$ shows \eqref{interpest_loc_s_q0plus_res}.
\end{proof}
\end{lemma}

\subsection{Global operator and estimates}

\begin{definition}(Main space-time interpolation operator)
We define the main space-time interpolation operator $\interph : C_b(\cup_n I_n^\circ; L^1(\Om{0})) \to V_h$ by, for $n = 1, \dots, N$,
\begin{equation}
\interph v|_{S_{1,n}} := \interp{{h,0}}{n} v|_{S_{1,n}} \qquad \interph v|_{S_{2,n}} := \interp{{h,G}}{n} v|_{S_{2,n}}
\label{def:interph}
\end{equation} 
\end{definition}

\begin{lemma}[Global space-time interpolation error estimates for $\Om{0} \times {(0, T]}$]
\label{lem:interpest_glob_st_om0T}
Let $\interph$ be defined by \eqref{def:interph} and $D_t$ by \eqref{def:Dt}. Then, for any function $v$ with sufficient spatial and temporal regularity, we have for $0 \leq s \leq q+1$ and $0 \leq r \leq p+1$ that
\begin{align}
\sum_{i=1}^2 \sum_{n=1}^N \int_{I_n} \| D_t^s ( v - \interph v) \|_{\Omt{i}}^2 \ud t & \lesssim k^{2(q+1-s)} E_{k, 0}^2(v) + h^{2(p+1)} E_{h, s}^2(v)
\label{interpest_glob_st_t_res} \\
\sum_{i=1}^2 \sum_{n=1}^N \int_{I_n} \| D_x^r ( v - \interph v) \|_{\Omt{i}}^2 \ud t & \lesssim k^{2(q+1)} E_{k, r}^2(v) + h^{2(p+1-r)} E_{h, 0}^2(v)
\label{interpest_glob_st_x_res}
\end{align}
where
\begin{align}
E_{k, r}^2(v) &= \sum_{i=1}^2 \sum_{n=1}^N \sum_{K \in \mathcal{T}_{j, S_{i,n}}} \| D_t^{q+1} D_x^r v \|_{\bar{K}_n}^2
\label{interpest_glob_st_Ek} \\
E_{h, s}^2(v) &= \sum_{i=1}^2 \sum_{n=1}^N \sum_{K \in \mathcal{T}_{j, S_{i,n}}} k_n \| D_x^{p+1} D_t^s v \|_{\mathcal{N}(K), I_n}^2
\label{interpest_glob_st_Eh}
\end{align}
\begin{proof}
Both estimates follow by applying Lemma~\ref{lem:interpest_loc_st_bKn}.
\end{proof}
\end{lemma}

\begin{lemma}[An interpolation error estimate in $\normb{\cdot}$] 
\label{lem:interpest}
Let $\normb{\cdot}$, $\interph$, and $D_t$ be defined by \eqref{def:bnorm}, \eqref{def:interph}, and \eqref{def:Dt}, respectively. Then, for any function $v$ with sufficient spatial and temporal regularity, we have that   
\begin{equation}
\begin{split}
\normb{v - \interph v}^2 \lesssim k^{2q +1} F_k^2(v) + h^{2p}F_h^2(v)
\end{split}
\label{interpest_res}   
\end{equation} 
where 
\begin{align}
\begin{split}
F_k^2(v) = & \; \sum_{i=1}^2\sum_{n=1}^N\sum_{K \in \mathcal{T}_{j, S_{i,n}}} \bigg(\| D_t^{q+1} v \|_{\bar{K}_n}^2 + \| D_t^{q+1} \nab v \|_{\bar{K}_n}^2 + \| D_t^{q+1}D_x^2 v \|_{\bar{K}_n}^2\bigg) \\
& + \sum_{i=1}^2\sum_{n=1}^N \| D_t^{q+1} v \|_{L^2(\bGn, L^\infty(I_n))}^2
\end{split}
\label{interpest_res_Fk} \\
F_h^2(v) = & \; \sum_{i=1}^2\sum_{n=1}^N\sum_{K \in \mathcal{T}_{j, S_{i,n}}} k_n \| D_x^{p+1} v \|_{\mathcal{N}(K), I_n}^2
\label{interpest_res_Fh}
\end{align}
\begin{proof} Letting $w = v - \interph v$, the left-hand side of \eqref{interpest_res} is 
\begin{equation}
\begin{split}
\normb{w}^2 = & \; \sum_{n=1}^N \underbrace{\int_{I_n} \norma{w}^2 \ud t}_\text{= I} + \sum_{n=1}^N \underbrace{\| |\bn^t|^{1/2} [w]\|_{\bGn}^2}_\text{= II} + \sum_{i=1}^2 \underbrace{\sum_{n=1}^{N-1} \|[w]_n \|_{\Om{{i,n}}}^2}_\text{= III} \\
& + \sum_{i=1}^2 \underbrace{\| w_N^- \|_{\Om{{i,N}}}^2}_\text{= IV} + \sum_{i=1}^2 \underbrace{ \| w_0^+ \|_{\Om{{i,0}}}^2}_\text{= V}
\end{split}
\label{interpest0}   
\end{equation} 
We consider the terms separately, starting with first. 
\begin{equation}
\begin{split}
\text{I} = & \; \sum_{i = 1}^2 \underbrace{\int_{I_n} \| \nab w \|_{\Omt{i}}^2 \ud t}_\text{= I.i} +  \underbrace{ \int_{I_n} \abm \|\langle \partial_{\bn^x} w \rangle \|_{-1/2,h,\Gamma(t)}^2 \ud t}_\text{= I.ii} \\
& +  \underbrace{ \int_{I_n} \abm \|[w] \|_{1/2,h,\Gamma(t)}^2 \ud t}_\text{= I.iii} +  \underbrace{ \int_{I_n} \|[\nab w]\|_{\Omt{O}}^2 \ud t}_\text{= I.iv}
\end{split}
\label{interpest0_1}   
\end{equation} 
Letting $w_j^n = v - \interp{{h, j}}{n}v$, we treat each term in (\ref{interpest0_1}) separately, starting with the first.
\begin{equation}
\text{I.i} \leq \int_{I_n} \sum_{K \in \mathcal{T}_{j, \Omt{i}}} \| \nab w_j^n \|_{K}^2 \ud t \leq \sum_{K \in \mathcal{T}_{j, S_{i,n}}} \| \nab  w_j^n \|_{\bar{K}_n}^2 
\label{interpest0_1_1}   
\end{equation} 
By using standard estimates, \eqref{gamk0ineqgamkg}, and \eqref{scatraineqGamK_warmup}, the second term is
\begin{equation}
\begin{split}
\text{I.ii} & \lesssim \int_{I_n} \sum_{i=1}^2 \sum_{K_j \in \mathcal{T}_{j, \Gamma(t)}} h_{K_j} \| (\nab w)_i \|_{\Gamma_{K_j}}^2 \ud t \\
& \lesssim \sum_{i=1}^2 \sum_{K \in \mathcal{T}_{j, \bGn}} \bigg(\| \nab w_j^n \|_{\bar{K}_n}^2 + h_K^2\| D_x^2 w_j^n \|_{\bar{K}_n}^2 \bigg)
\end{split}
\label{interpest0_1_2}   
\end{equation} 
For the third term we use the same standard estimates and again \eqref{gamk0ineqgamkg}, thus 
\begin{equation}
\text{I.iii} \lesssim \int_{I_n} \abm \sum_{i=1}^2 \sum_{K_j \in \mathcal{T}_{j, \Gamma(t)}} h_{K_0}^{-1} \| w_i \|_{\Gamma_{K_j}}^2 \ud t \leq  h_{\min}^{-1} \sum_{i=1}^2 \| (w_j^n)_i \|_{\bGn}^2
\label{interpest0_1_3}   
\end{equation} 
The fourth term is
\begin{equation}
\text{I.iv} \lesssim \int_{I_n} \sum_{i=1}^2\|(\nab w)_i \|_{\Omt{O}}^2 \ud t \leq \sum_{i=1}^2 \sum_{K \in \mathcal{T}_{j, \bGn}} \| \nab w_j^n \|_{\bar{K}_n}^2
\label{interpest0_1_4}   
\end{equation} 
We are done with the separate treatments of all the terms in \eqref{interpest0_1} and move on to the second term in \eqref{interpest0}. For this term, using that $|\bn^t| \leq |\mu|$ and \eqref{HHnormineqs} results in a factor that is I.iii which may simply be estimated by \eqref{interpest0_1_3}, thus
\begin{equation}
\text{II} \leq | \mu |_{(0,T]} h \int_{I_n} \abm \|[w]\|_{1/2,h,\Gamma(t)}^2 \ud t \lesssim | \mu |_{(0,T]} h h_{\min}^{-1}\sum_{i=1}^2 \| (w_j^n)_i \|_{\bGn}^2
\label{interpest0_2}   
\end{equation} 
Combining the third, fourth and fifth term in \eqref{interpest0}, we have
\begin{equation}
\begin{split}
\text{III + IV + V} & \lesssim \sum_{n=1}^{N} \bigg( \| w_n^- \|_{\Om{{i,n}}}^2 + \| w_{n-1}^+ \|_{\Om{{i,n-1}}}^2 \bigg) \\
& \leq \sum_{n=1}^{N} \bigg( \sum_{K \in \mathcal{T}_{j, \Om{{i,n}}}} \| (w_j^n)_n^- \|_{K}^2 + \sum_{K \in \mathcal{T}_{j, \Om{{i,n-1}}}} \| (w_j^n)_{n-1}^+ \|_{K}^2 \bigg)
\end{split}
\label{interpest0_3}
\end{equation}
The separate treatments of all the terms in \eqref{interpest0} are done. The obtained estimates give
\begin{equation}
\begin{split}
\normb{w}^2 \lesssim & \; \sum_{n=1}^N \sum_{i=1}^2 \bigg( \sum_{K \in \mathcal{T}_{j, S_{i,n}}} \underbrace{\| \nab  w_j^n \|_{\bar{K}_n}^2}_\text{= A} +  \sum_{K \in \mathcal{T}_{j, \bGn}} \underbrace{h_K^2\| D_x^2 w_j^n \|_{\bar{K}_n}^2}_\text{= B} \\
& +  \underbrace{h_{\min}^{-1} \| (w_j^n)_i \|_{\bGn}^2}_\text{= C} + \sum_{K \in \mathcal{T}_{j, \Om{{i,n}}}} \underbrace{\| (w_j^n)_n^- \|_{K}^2}_\text{= D} + \sum_{K \in \mathcal{T}_{j, \Om{{i,n-1}}}} \underbrace{\| (w_j^n)_{n-1}^+ \|_{K}^2}_\text{= E} \bigg) \\
\end{split}
\label{interpest_mid}   
\end{equation} 
We proceed by considering the five different types of terms separately. For term A we use Lemma~\ref{lem:interpest_loc_st_bKn} with $r = 1$:
\begin{equation}
\text{A} = \|\nab(v - \interp{{h,j}}{n}v)\|_{\bar{K}_n}^2 \lesssim k_n^{2(q + 1)} \| D_t^{q+1} \nab v \|_{\bar{K}_n}^2 + h^{2p} k_n \| D_x^{p+1} v \|_{\mathcal{N}(K), I_n}^2
\label{interpest_mid_A}
\end{equation}
For term B we apply Lemma~\ref{lem:interpest_loc_st_bKn} with $r = 2$:
\begin{equation}
\text{B} = h_{K}^{2} \| D_x^2(v - \interp{{h,j}}{n}v) \|_{\bar{K}_n}^2 \lesssim k_n^{2(q+1)}\| D_t^{q+1}D_x^2 v \|_{\bar{K}_n}^2 + h^{2p} k_n \| D_x^{p+1} v \|_{\mathcal{N}(K), I_n}^2
\label{interpest_mid_B}
\end{equation}
For term C we use Lemma~\ref{lem:interpest_loc_st_bGn} and \eqref{quasiuniformity_st}:
\begin{equation}
\begin{split}
\text{C} &= h_{\min}^{-1} \| (v - \interp{{h,j}}{n} v)_i \|_{\bGn}^2 \\ 
&\lesssim k_n^{2q+1} \| D_t^{q+1} v \|_{L^2(\bGn, L^\infty(I_n))}^2 + h^{2p}  \sum_{K \in \mathcal{T}_{j, \bGn}} k_n \| D_x^{p+1}v \|_{\mathcal{N}(K), I_n}^2
\end{split}
\label{interpest_mid_C}
\end{equation}
By applying Lemma~\ref{lem:interpest_loc_s} to term D and using \eqref{quasiuniformity_st}, we get
\begin{equation}
\begin{split}
\text{D} = \|(v - \interp{{h,j}}{n}v)_n^- \|_{K}^2 \lesssim h^{2p} k_n \| D_x^{p+1}v \|_{\mathcal{N}(K), I_n}^2
\end{split}
\label{interpest_mid_D}
\end{equation}
Using Lemma~\ref{lem:interpest_loc_s} and \eqref{quasiuniformity_st} for term E, we get
\begin{equation}
\text{E} = \|(v - \interp{{h,j}}{n}v)_{n-1}^+ \|_{K}^2 \lesssim k_n^{2q + 1} \| D_t^{q+1} v \|_{\bar{K}_n}^2 + h^{2p} k_n \| D_x^{p+1}v \|_{\mathcal{N}(K), I_n}^2
\label{interpest_mid_Eqcomb}
\end{equation}
Using these local estimates in \eqref{interpest_mid} gives \eqref{interpest_res}. 
\end{proof}
\end{lemma}
\noindent By applying Lemma~\ref{lem:interpest_loc_s} and using \eqref{quasiuniformity_st}, we get the estimate:
\begin{corollary}[A global spatial interpolation error estimate for temporal endpoints]
\label{cor:interpest_glob_s}
Let $\interph$ and $F_h$ be defined by \eqref{def:interph} and \eqref{interpest_res_Fh}, respectively. Then, for any function $v$ with sufficient spatial and temporal regularity, we have that 
\begin{equation}
\sum_{n=1}^N \| (v - \interph v)_n^- \|_{\Om{0}}^2 \lesssim h^{2p}F_h^2(v)
\label{interpest_glob_s_res}
\end{equation}
\end{corollary}

\bibliographystyle{IEEEtran}
\bibliography{bibliography}

\end{document}

%% file: figures/problemdomains.pdf_tex
\begingroup%
  \makeatletter%
  \providecommand\color[2][]{%
    \errmessage{(Inkscape) Color is used for the text in Inkscape, but the package 'color.sty' is not loaded}%
    \renewcommand\color[2][]{}%
  }%
  \providecommand\transparent[1]{%
    \errmessage{(Inkscape) Transparency is used (non-zero) for the text in Inkscape, but the package 'transparent.sty' is not loaded}%
    \renewcommand\transparent[1]{}%
  }%
  \providecommand\rotatebox[2]{#2}%
  \newcommand*\fsize{\dimexpr\f@size pt\relax}%
  \newcommand*\lineheight[1]{\fontsize{\fsize}{#1\fsize}\selectfont}%
  \ifx\svgwidth\undefined%
    \setlength{\unitlength}{1348.68937852bp}%
    \ifx\svgscale\undefined%
      \relax%
    \else%
      \setlength{\unitlength}{\unitlength * \real{\svgscale}}%
    \fi%
  \else%
    \setlength{\unitlength}{\svgwidth}%
  \fi%
  \global\let\svgwidth\undefined%
  \global\let\svgscale\undefined%
  \makeatother%
  \begin{picture}(1,0.62242716)%
    \lineheight{1}%
    \setlength\tabcolsep{0pt}%
    \put(0.78294052,0.16163747){\color[rgb]{0,0,0}\makebox(0,0)[lt]{\lineheight{0}\smash{\begin{tabular}[t]{l}\large : $\Omega_1$\end{tabular}}}}%
    \put(0,0){\includegraphics[width=\unitlength,page=1]{problemdomains.pdf}}%
    \put(0.78294052,0.08378411){\color[rgb]{0,0,0}\makebox(0,0)[lt]{\lineheight{0}\smash{\begin{tabular}[t]{l}\large : $\Omega_2$\end{tabular}}}}%
    \put(0.78294052,0.00593075){\color[rgb]{0,0,0}\makebox(0,0)[lt]{\lineheight{0}\smash{\begin{tabular}[t]{l}\large : $\Gamma$\end{tabular}}}}%
    \put(0,0){\includegraphics[width=\unitlength,page=2]{problemdomains.pdf}}%
    \put(0.5003338,0.07073026){\color[rgb]{0,0,0}\makebox(0,0)[lt]{\lineheight{0}\smash{\begin{tabular}[t]{l}\large $\Omega_0$\end{tabular}}}}%
    \put(0.37285761,0.0923501){\color[rgb]{0,0,0}\makebox(0,0)[lt]{\lineheight{0}\smash{\begin{tabular}[t]{l}\large $\mu$\end{tabular}}}}%
    \put(0.31234314,0.29020725){\color[rgb]{0,0,0}\makebox(0,0)[lt]{\lineheight{0}\smash{\begin{tabular}[t]{l}\large $G$\end{tabular}}}}%
    \put(0,0){\includegraphics[width=\unitlength,page=3]{problemdomains.pdf}}%
  \end{picture}%
\endgroup%

%% file: figures/spacetimeslabs_cG1mesh.pdf_tex
\begingroup%
  \makeatletter%
  \providecommand\color[2][]{%
    \errmessage{(Inkscape) Color is used for the text in Inkscape, but the package 'color.sty' is not loaded}%
    \renewcommand\color[2][]{}%
  }%
  \providecommand\transparent[1]{%
    \errmessage{(Inkscape) Transparency is used (non-zero) for the text in Inkscape, but the package 'transparent.sty' is not loaded}%
    \renewcommand\transparent[1]{}%
  }%
  \providecommand\rotatebox[2]{#2}%
  \newcommand*\fsize{\dimexpr\f@size pt\relax}%
  \newcommand*\lineheight[1]{\fontsize{\fsize}{#1\fsize}\selectfont}%
  \ifx\svgwidth\undefined%
    \setlength{\unitlength}{2069.29439705bp}%
    \ifx\svgscale\undefined%
      \relax%
    \else%
      \setlength{\unitlength}{\unitlength * \real{\svgscale}}%
    \fi%
  \else%
    \setlength{\unitlength}{\svgwidth}%
  \fi%
  \global\let\svgwidth\undefined%
  \global\let\svgscale\undefined%
  \makeatother%
  \begin{picture}(1,0.58581438)%
    \lineheight{1}%
    \setlength\tabcolsep{0pt}%
    \put(0.80868287,0.17397535){\color[rgb]{0,0,0}\makebox(0,0)[lt]{\lineheight{0}\smash{\begin{tabular}[t]{l}$S_{1,n-1}$\end{tabular}}}}%
    \put(0,0){\includegraphics[width=\unitlength,page=1]{spacetimeslabs_cG1mesh.pdf}}%
    \put(0.16278541,0.08884568){\color[rgb]{0,0,0}\makebox(0,0)[lt]{\lineheight{0}\smash{\begin{tabular}[t]{l}$x_2$\end{tabular}}}}%
    \put(0,0){\includegraphics[width=\unitlength,page=2]{spacetimeslabs_cG1mesh.pdf}}%
    \put(0.18791888,0.00217465){\color[rgb]{0,0,0}\makebox(0,0)[lt]{\lineheight{0}\smash{\begin{tabular}[t]{l}$x_1$\end{tabular}}}}%
    \put(-0.00076548,0.19035702){\color[rgb]{0,0,0}\makebox(0,0)[lt]{\lineheight{0}\smash{\begin{tabular}[t]{l}$t$\end{tabular}}}}%
    \put(0,0){\includegraphics[width=\unitlength,page=3]{spacetimeslabs_cG1mesh.pdf}}%
    \put(0.48322076,0.16962897){\color[rgb]{0,0,0}\makebox(0,0)[lt]{\lineheight{0}\smash{\begin{tabular}[t]{l}$S_{2,n-1}$\end{tabular}}}}%
    \put(0.14772529,0.03866927){\color[rgb]{0,0,0}\makebox(0,0)[lt]{\lineheight{0}\smash{\begin{tabular}[t]{l}$t_{n-2}$\end{tabular}}}}%
    \put(0,0){\includegraphics[width=\unitlength,page=4]{spacetimeslabs_cG1mesh.pdf}}%
    \put(0.80868287,0.40593848){\color[rgb]{0,0,0}\makebox(0,0)[lt]{\lineheight{0}\smash{\begin{tabular}[t]{l}$S_{1,n}$\end{tabular}}}}%
    \put(0,0){\includegraphics[width=\unitlength,page=5]{spacetimeslabs_cG1mesh.pdf}}%
    \put(0.17607703,0.50251425){\color[rgb]{0,0,0}\makebox(0,0)[lt]{\lineheight{0}\smash{\begin{tabular}[t]{l}$t_n$\end{tabular}}}}%
    \put(0,0){\includegraphics[width=\unitlength,page=6]{spacetimeslabs_cG1mesh.pdf}}%
    \put(0.63544656,0.40884094){\color[rgb]{0,0,0}\makebox(0,0)[lt]{\lineheight{0}\smash{\begin{tabular}[t]{l}$S_{2,n}$\end{tabular}}}}%
    \put(0.14772529,0.27063239){\color[rgb]{0,0,0}\makebox(0,0)[lt]{\lineheight{0}\smash{\begin{tabular}[t]{l}$t_{n-1}$\end{tabular}}}}%
    \put(0,0){\includegraphics[width=\unitlength,page=7]{spacetimeslabs_cG1mesh.pdf}}%
  \end{picture}%
\endgroup%

%% file: figures/rumtidsmesh1dmup.pdf_tex
\begingroup%
  \makeatletter%
  \providecommand\color[2][]{%
    \errmessage{(Inkscape) Color is used for the text in Inkscape, but the package 'color.sty' is not loaded}%
    \renewcommand\color[2][]{}%
  }%
  \providecommand\transparent[1]{%
    \errmessage{(Inkscape) Transparency is used (non-zero) for the text in Inkscape, but the package 'transparent.sty' is not loaded}%
    \renewcommand\transparent[1]{}%
  }%
  \providecommand\rotatebox[2]{#2}%
  \ifx\svgwidth\undefined%
    \setlength{\unitlength}{877.1702368bp}%
    \ifx\svgscale\undefined%
      \relax%
    \else%
      \setlength{\unitlength}{\unitlength * \real{\svgscale}}%
    \fi%
  \else%
    \setlength{\unitlength}{\svgwidth}%
  \fi%
  \global\let\svgwidth\undefined%
  \global\let\svgscale\undefined%
  \makeatother%
  \begin{picture}(1,0.48558324)%
    \put(0,0){\includegraphics[width=\unitlength]{rumtidsmesh1dmup.pdf}}%
    \put(-0.00179468,0.07898823){\color[rgb]{0,0,0}\makebox(0,0)[lt]{\begin{minipage}{0.13666193\unitlength}\raggedright $t_{n-1}$\end{minipage}}}%
    \put(0.30642601,0.00250094){\color[rgb]{0,0,0}\makebox(0,0)[lb]{\smash{$x$}}}%
    \put(0.07889132,0.46107083){\color[rgb]{0,0,0}\makebox(0,0)[lb]{\smash{$t$}}}%
    \put(-0.00179468,0.35830354){\color[rgb]{0,0,0}\makebox(0,0)[lt]{\begin{minipage}{0.13666193\unitlength}\raggedright $t_n$\end{minipage}}}%
  \end{picture}%
\endgroup%

%% file: figures/spacetimevector.pdf_tex
\begingroup%
  \makeatletter%
  \providecommand\color[2][]{%
    \errmessage{(Inkscape) Color is used for the text in Inkscape, but the package 'color.sty' is not loaded}%
    \renewcommand\color[2][]{}%
  }%
  \providecommand\transparent[1]{%
    \errmessage{(Inkscape) Transparency is used (non-zero) for the text in Inkscape, but the package 'transparent.sty' is not loaded}%
    \renewcommand\transparent[1]{}%
  }%
  \providecommand\rotatebox[2]{#2}%
  \newcommand*\fsize{\dimexpr\f@size pt\relax}%
  \newcommand*\lineheight[1]{\fontsize{\fsize}{#1\fsize}\selectfont}%
  \ifx\svgwidth\undefined%
    \setlength{\unitlength}{1449.70919861bp}%
    \ifx\svgscale\undefined%
      \relax%
    \else%
      \setlength{\unitlength}{\unitlength * \real{\svgscale}}%
    \fi%
  \else%
    \setlength{\unitlength}{\svgwidth}%
  \fi%
  \global\let\svgwidth\undefined%
  \global\let\svgscale\undefined%
  \makeatother%
  \begin{picture}(1,0.49412723)%
    \lineheight{1}%
    \setlength\tabcolsep{0pt}%
    \put(0,0){\includegraphics[width=\unitlength,page=1]{spacetimevector.pdf}}%
    \put(0.25564722,0.12681707){\color[rgb]{0,0,0}\makebox(0,0)[lt]{\lineheight{0}\smash{\begin{tabular}[t]{l}$x_2$\end{tabular}}}}%
    \put(0,0){\includegraphics[width=\unitlength,page=2]{spacetimevector.pdf}}%
    \put(0.29152237,0.00310407){\color[rgb]{0,0,0}\makebox(0,0)[lt]{\lineheight{0}\smash{\begin{tabular}[t]{l}$x_1$\end{tabular}}}}%
    \put(0.33322701,0.15361121){\color[rgb]{0,0,0}\makebox(0,0)[lt]{\lineheight{0}\smash{\begin{tabular}[t]{l}$n_2$\end{tabular}}}}%
    \put(0.38595756,0.26915209){\color[rgb]{0,0,0}\makebox(0,0)[lt]{\lineheight{0}\smash{\begin{tabular}[t]{l}$\bar{n}_2$\end{tabular}}}}%
    \put(0.69464361,0.27846418){\color[rgb]{0,0,0}\makebox(0,0)[lt]{\lineheight{0}\smash{\begin{tabular}[t]{l}\large $\bar{\Gamma}_n$\end{tabular}}}}%
    \put(0.62894987,0.20179454){\color[rgb]{0,0,0}\makebox(0,0)[lt]{\lineheight{0}\smash{\begin{tabular}[t]{l}\large $\Omega_2(t^*)$\end{tabular}}}}%
    \put(0,0){\includegraphics[width=\unitlength,page=3]{spacetimevector.pdf}}%
    \put(-0.00109263,0.05251805){\color[rgb]{0,0,0}\makebox(0,0)[lt]{\lineheight{0}\smash{\begin{tabular}[t]{l}$t_{n-1}$\end{tabular}}}}%
    \put(0,0){\includegraphics[width=\unitlength,page=4]{spacetimevector.pdf}}%
    \put(0.02777992,0.40269673){\color[rgb]{0,0,0}\makebox(0,0)[lt]{\lineheight{0}\smash{\begin{tabular}[t]{l}$t_n$\end{tabular}}}}%
    \put(0,0){\includegraphics[width=\unitlength,page=5]{spacetimevector.pdf}}%
    \put(0.02891338,0.19737468){\color[rgb]{0,0,0}\makebox(0,0)[lt]{\lineheight{0}\smash{\begin{tabular}[t]{l}$t^*$\end{tabular}}}}%
    \put(0.03502276,0.45339655){\color[rgb]{0,0,0}\makebox(0,0)[lt]{\lineheight{0}\smash{\begin{tabular}[t]{l}$t$\end{tabular}}}}%
    \put(0,0){\includegraphics[width=\unitlength,page=6]{spacetimevector.pdf}}%
  \end{picture}%
\endgroup%

%% file: figures/fefuncspace.pdf_tex
\begingroup%
  \makeatletter%
  \providecommand\color[2][]{%
    \errmessage{(Inkscape) Color is used for the text in Inkscape, but the package 'color.sty' is not loaded}%
    \renewcommand\color[2][]{}%
  }%
  \providecommand\transparent[1]{%
    \errmessage{(Inkscape) Transparency is used (non-zero) for the text in Inkscape, but the package 'transparent.sty' is not loaded}%
    \renewcommand\transparent[1]{}%
  }%
  \providecommand\rotatebox[2]{#2}%
  \ifx\svgwidth\undefined%
    \setlength{\unitlength}{861.750576bp}%
    \ifx\svgscale\undefined%
      \relax%
    \else%
      \setlength{\unitlength}{\unitlength * \real{\svgscale}}%
    \fi%
  \else%
    \setlength{\unitlength}{\svgwidth}%
  \fi%
  \global\let\svgwidth\undefined%
  \global\let\svgscale\undefined%
  \makeatother%
  \begin{picture}(1,0.49427199)%
    \put(0,0){\includegraphics[width=\unitlength]{fefuncspace.pdf}}%
    \put(0.2940156,0.00254569){\color[rgb]{0,0,0}\makebox(0,0)[lb]{\smash{$x$}}}%
    \put(0.06240954,0.46932096){\color[rgb]{0,0,0}\makebox(0,0)[lb]{\smash{$v(x,t)$}}}%
    \put(-0.0018268,0.06169315){\color[rgb]{0,0,0}\makebox(0,0)[lt]{\begin{minipage}{0.13910728\unitlength}\raggedright $0$\end{minipage}}}%
  \end{picture}%
\endgroup%

%% file: figures/auxfig_lemma_disctemp_invest.pdf_tex
\begingroup%
  \makeatletter%
  \providecommand\color[2][]{%
    \errmessage{(Inkscape) Color is used for the text in Inkscape, but the package 'color.sty' is not loaded}%
    \renewcommand\color[2][]{}%
  }%
  \providecommand\transparent[1]{%
    \errmessage{(Inkscape) Transparency is used (non-zero) for the text in Inkscape, but the package 'transparent.sty' is not loaded}%
    \renewcommand\transparent[1]{}%
  }%
  \providecommand\rotatebox[2]{#2}%
  \newcommand*\fsize{\dimexpr\f@size pt\relax}%
  \newcommand*\lineheight[1]{\fontsize{\fsize}{#1\fsize}\selectfont}%
  \ifx\svgwidth\undefined%
    \setlength{\unitlength}{820.25999645bp}%
    \ifx\svgscale\undefined%
      \relax%
    \else%
      \setlength{\unitlength}{\unitlength * \real{\svgscale}}%
    \fi%
  \else%
    \setlength{\unitlength}{\svgwidth}%
  \fi%
  \global\let\svgwidth\undefined%
  \global\let\svgscale\undefined%
  \makeatother%
  \begin{picture}(1,0.43639877)%
    \lineheight{1}%
    \setlength\tabcolsep{0pt}%
    \put(0,0){\includegraphics[width=\unitlength,page=1]{auxfig_lemma_disctemp_invest.pdf}}%
    \put(0.30179614,0.00452016){\color[rgb]{0,0,0}\makebox(0,0)[lt]{\lineheight{0}\smash{\begin{tabular}[t]{l}$x$\end{tabular}}}}%
    \put(0.05451299,0.40311245){\color[rgb]{0,0,0}\makebox(0,0)[lt]{\lineheight{0}\smash{\begin{tabular}[t]{l}$t$\end{tabular}}}}%
    \put(0,0){\includegraphics[width=\unitlength,page=2]{auxfig_lemma_disctemp_invest.pdf}}%
    \put(0.03134184,0.34622946){\color[rgb]{0,0,0}\makebox(0,0)[lt]{\lineheight{0}\smash{\begin{tabular}[t]{l}$t_n$\end{tabular}}}}%
    \put(-0.00144832,0.06587381){\color[rgb]{0,0,0}\makebox(0,0)[lt]{\lineheight{0}\smash{\begin{tabular}[t]{l}$t_{n-1}$\end{tabular}}}}%
    \put(0.77176911,0.18373163){\color[rgb]{0,0,0}\makebox(0,0)[lt]{\lineheight{0}\smash{\begin{tabular}[t]{l}$S_{1,n}$\end{tabular}}}}%
    \put(0.41151745,0.18373163){\color[rgb]{0,0,0}\makebox(0,0)[lt]{\lineheight{0}\smash{\begin{tabular}[t]{l}$S_{2,n}$\end{tabular}}}}%
    \put(0.57975682,0.12338491){\color[rgb]{0,0,0}\makebox(0,0)[lt]{\lineheight{0}\smash{\begin{tabular}[t]{l}$\bar{s}$\end{tabular}}}}%
    \put(0.56146993,0.02280705){\color[rgb]{0,0,0}\makebox(0,0)[lt]{\lineheight{0}\smash{\begin{tabular}[t]{l}$(\hat{s}_{1,n-1}, t_{n-1})$\end{tabular}}}}%
    \put(0.35848551,0.02280705){\color[rgb]{0,0,0}\makebox(0,0)[lt]{\lineheight{0}\smash{\begin{tabular}[t]{l}$(\hat{s}_{2,n-1}, t_{n-1})$\end{tabular}}}}%
    \put(0.66021911,0.36660049){\color[rgb]{0,0,0}\makebox(0,0)[lt]{\lineheight{0}\smash{\begin{tabular}[t]{l}$(\hat{s}_{2,n}, t_n)$\end{tabular}}}}%
    \put(0.45174863,0.36660049){\color[rgb]{0,0,0}\makebox(0,0)[lt]{\lineheight{0}\smash{\begin{tabular}[t]{l}$(\hat{s}_{1,n}, t_n)$\end{tabular}}}}%
    \put(0,0){\includegraphics[width=\unitlength,page=3]{auxfig_lemma_disctemp_invest.pdf}}%
  \end{picture}%
\endgroup%

%% file: figures/auxfig_lemma_disctemp_invineq.pdf_tex
\begingroup%
  \makeatletter%
  \providecommand\color[2][]{%
    \errmessage{(Inkscape) Color is used for the text in Inkscape, but the package 'color.sty' is not loaded}%
    \renewcommand\color[2][]{}%
  }%
  \providecommand\transparent[1]{%
    \errmessage{(Inkscape) Transparency is used (non-zero) for the text in Inkscape, but the package 'transparent.sty' is not loaded}%
    \renewcommand\transparent[1]{}%
  }%
  \providecommand\rotatebox[2]{#2}%
  \newcommand*\fsize{\dimexpr\f@size pt\relax}%
  \newcommand*\lineheight[1]{\fontsize{\fsize}{#1\fsize}\selectfont}%
  \ifx\svgwidth\undefined%
    \setlength{\unitlength}{927.84506085bp}%
    \ifx\svgscale\undefined%
      \relax%
    \else%
      \setlength{\unitlength}{\unitlength * \real{\svgscale}}%
    \fi%
  \else%
    \setlength{\unitlength}{\svgwidth}%
  \fi%
  \global\let\svgwidth\undefined%
  \global\let\svgscale\undefined%
  \makeatother%
  \begin{picture}(1,0.38579766)%
    \lineheight{1}%
    \setlength\tabcolsep{0pt}%
    \put(0,0){\includegraphics[width=\unitlength,page=1]{auxfig_lemma_disctemp_invineq.pdf}}%
    \put(0.26841906,0.00399604){\color[rgb]{0,0,0}\makebox(0,0)[lt]{\lineheight{0}\smash{\begin{tabular}[t]{l}$x$\end{tabular}}}}%
    \put(0.04819212,0.35637094){\color[rgb]{0,0,0}\makebox(0,0)[lt]{\lineheight{0}\smash{\begin{tabular}[t]{l}$t$\end{tabular}}}}%
    \put(0,0){\includegraphics[width=\unitlength,page=2]{auxfig_lemma_disctemp_invineq.pdf}}%
    \put(0.0277077,0.30608362){\color[rgb]{0,0,0}\makebox(0,0)[lt]{\lineheight{0}\smash{\begin{tabular}[t]{l}$t_n$\end{tabular}}}}%
    \put(0,0){\includegraphics[width=\unitlength,page=3]{auxfig_lemma_disctemp_invineq.pdf}}%
    \put(-0.00128039,0.05823564){\color[rgb]{0,0,0}\makebox(0,0)[lt]{\lineheight{0}\smash{\begin{tabular}[t]{l}$t_{n-1}$\end{tabular}}}}%
    \put(0,0){\includegraphics[width=\unitlength,page=4]{auxfig_lemma_disctemp_invineq.pdf}}%
    \put(0.16694897,0.01369594){\color[rgb]{0,0,0}\makebox(0,0)[lt]{\lineheight{0}\smash{\begin{tabular}[t]{l}$\Omega_{1,n-1}^{n}$\end{tabular}}}}%
    \put(0.14878702,0.17859416){\color[rgb]{0,0,0}\makebox(0,0)[lt]{\lineheight{0}\smash{\begin{tabular}[t]{l}I\end{tabular}}}}%
    \put(0.34763488,0.17859416){\color[rgb]{0,0,0}\makebox(0,0)[lt]{\lineheight{0}\smash{\begin{tabular}[t]{l}II\end{tabular}}}}%
    \put(0.44463383,0.16242767){\color[rgb]{0,0,0}\makebox(0,0)[lt]{\lineheight{0}\smash{\begin{tabular}[t]{l}III.3\end{tabular}}}}%
    \put(0.49313331,0.24164348){\color[rgb]{0,0,0}\makebox(0,0)[lt]{\lineheight{0}\smash{\begin{tabular}[t]{l}III.2\end{tabular}}}}%
    \put(0.52869959,0.09937835){\color[rgb]{0,0,0}\makebox(0,0)[lt]{\lineheight{0}\smash{\begin{tabular}[t]{l}III.1\end{tabular}}}}%
    \put(0.33631833,0.01369594){\color[rgb]{0,0,0}\makebox(0,0)[lt]{\lineheight{0}\smash{\begin{tabular}[t]{l}$\Omega_{2,n-1}$\end{tabular}}}}%
    \put(0.51414975,0.01369594){\color[rgb]{0,0,0}\makebox(0,0)[lt]{\lineheight{0}\smash{\begin{tabular}[t]{l}$\Omega_{1,n-1}^{\Gamma}$\end{tabular}}}}%
    \put(0.72108087,0.01369594){\color[rgb]{0,0,0}\makebox(0,0)[lt]{\lineheight{0}\smash{\begin{tabular}[t]{l}$\Omega_{1,n-1}^{n}$\end{tabular}}}}%
    \put(0,0){\includegraphics[width=\unitlength,page=5]{auxfig_lemma_disctemp_invineq.pdf}}%
  \end{picture}%
\endgroup%